\documentclass{unbranded}
\pdfoutput=1

\OneAndAHalfSpacedXII

\usepackage[utf8]{inputenc}
\usepackage{algorithm}
\usepackage{algpseudocode}
\usepackage{amsmath,amssymb}
\usepackage{array}
\usepackage{calc}
\usepackage{enumitem}
\usepackage[T1]{fontenc}
\usepackage[margin=1.00in]{geometry}
\usepackage{hhline}
\usepackage{makecell}
\usepackage{multirow,multicol}
\usepackage[short,nocomma]{optidef}
\usepackage{subcaption}
\usepackage[title]{appendix}  % outcomment this for submission
\usepackage{hyperref}

\hypersetup{colorlinks=true,urlcolor=blue,linkcolor=blue,citecolor=red}

\DeclareMathOperator{\prob}{Pr}
\def\tigerdam{Tiger Dam\texttrademark}
\def\tigerdams{Tiger Dams\texttrademark}

\def\bs{\boldsymbol}
\newcolumntype{L}[1]{>{\raggedright\let\newline\\\arraybackslash\hspace{0pt}}m{#1}}
\newcolumntype{C}[1]{>{\centering\let\newline\\\arraybackslash\hspace{0pt}}m{#1}}
\newcolumntype{R}[1]{>{\raggedleft\let\newline\\\arraybackslash\hspace{0pt}}m{#1}}
\DeclareMathOperator*{\abssim}{AbsSim}
\DeclareMathOperator*{\relsim}{RelSim}
\DeclareMathOperator{\diag}{diag}
\let\epsswitch\undefined  % set to false

\usepackage[normalem]{ulem}
\usepackage{censor}
\newlength\nextcharwidth
\makeatletter
\renewcommand\@cenword[1]{%
  \setlength{\nextcharwidth}{\widthof{#1}}%
  \censorrule{\nextcharwidth}%
  \kern -\nextcharwidth%
  #1}
\makeatother

\definecolor{response}{RGB}{0, 0, 255}
\definecolor{strikeout}{RGB}{255, 0, 0}

\usepackage{natbib}
 \bibpunct[, ]{(}{)}{,}{a}{}{,}%
 \def\bibfont{\small}%
\TheoremsNumberedThrough     % Preferred (Theorem 1, Lemma 1, Theorem 2)
\EquationsNumberedThrough    % Default: (1), (2), ...
\begin{document}
%%%%%%%%%%%%%%%%

%\setcounter{page}{1}
\DOI{10.1287/ijoc.2023.0125}%

\RUNAUTHOR{Austgen et al.}
\RUNTITLE{Comparisons of Two-stage Models for Flood Mitigation of Electrical Substations}
\TITLE{Comparisons of Two-stage Models for Flood Mitigation of Electrical Substations}

\ARTICLEAUTHORS{%
    \AUTHOR{Brent Austgen,
            Erhan Kutanoglu,
            John J. Hasenbein}
    \AFF{The University of Texas at Austin, Operations Research and Industrial Engineering Program}
    \AUTHOR{Surya Santoso}
    \AFF{The University of Texas at Austin, Chandra Department of Electrical and Computer Engineering}
}

\ABSTRACT{%
We compare stochastic programming and robust optimization decision models for informing the deployment of ad hoc flood mitigation measures to protect electrical substations prior to an imminent and uncertain hurricane. In our models, the first stage captures the deployment of a fixed quantity of flood mitigation resources, and the second stage captures the operation of a potentially degraded power grid with the primary goal of minimizing load shed. To model grid operation, we introduce adaptations of the DC and LPAC power flow approximation models that feature relatively complete recourse by way of an indicator variable. We apply our models to a pair of geographically realistic flooding case studies, one based on Hurricane Harvey and the other on Tropical Storm Imelda. We investigate the effect of the mitigation budget, the choice of power flow model, and the uncertainty perspective on the optimal mitigation strategy. Our results indicate the mitigation budget and uncertainty perspective are impactful whereas choosing between the DC and LPAC power flow models is of little to no consequence. To validate our models, we assess the performance of the mitigation solutions they prescribe in an AC power flow model.
}%


\maketitle

% Text of your paper here
\section{Introduction} \label{section:intro}

The power grid underpins or is codependent with many lifeline infrastructure systems including water, natural gas, transportation, and communication. Power grid resilience is obviously crucial to humankind's well-being; however, it is also vulnerable especially to extreme weather. According to information collected by the United States Department of Energy via form DOE-417, the number of widespread outages caused by severe weather and natural disasters has increased from around 40 annually in the early 2000's to about 80 by 2010 to 100 or more in recent years \citep{ISER2023}. Many of the worst outages and economic losses are due to tropical cyclones \citep{BillionDollarDisasters2022}, and many climate models project the frequency and intensity of the most extreme tropical cyclones (\textit{i.e.}, Category 4 and Category 5 storms) to increase globally \citep{Webster2005,Knutson2020}, thus motivating the continued study and pursuit of power grid resilience.

In decision models, power grid physics and operation are captured by a power flow (PF) model. The exact model for AC power flow (ACPF) is nonconvex and thus computationally difficult in many decision making applications. However, a variety of tractable surrogate models may be formed by approximating, relaxing, or restricting the exact model \citep{Molzahn2019}. In deciding which PF model to use for an application, the trade-off between a model's fidelity and computational ease is a major consideration.

In our application, we consider the mitigation of flood-induced component outages via deployment of temporary flood barriers like \tigerdams~in the hours leading up to an imminent hurricane's landfall. \tigerdams~are rapidly deployable mitigation resources that may be connected end-to-end to lengthen a barrier and stacked to heighten it. The discrete mitigation decisions and flooding uncertainty induce combinatorial effects on the grid topology which require the conditional enforcement of certain physical laws (\textit{e.g.}, only enforce Ohm's Law on closed circuits).
Due to these inherent complexities and the time-sensitive nature of the decision making in our application, we only consider incorporating simple linear and convex quadratic PF approximations.
Even among these simple PF models, there is considerable variation in fidelity and complexity.

In this paper, we extend the past work to deliver the following contributions:
\begin{enumerate}
    \item We formulate alternatives to the two-stage stochastic programming (SP) model with a recourse problem based on DC power flow that we proposed for substation flood mitigation in \citet{Austgen2023}. Specifically, we consider an analogous two-stage robust optimization (RO) formulation and recourse problems based on three dynamic-topology variants of the linear programming AC (LPAC) power flow model \citep{Coffrin2014}.
    \item We extend our recourse problem from \citet{Austgen2023} by incorporating a binary decision that, if exercised, allows a trivial recourse problem solution in which no power is generated and no loads are satisfied. We present this ``infeasibility indicator variable'' as a coarse alternative to relatively complicated discrete controls like transmission line switches and generator commitments that may similarly ensure power flow feasibility.
    \item In a computational study, we apply our proposed two-stage models to two geographically realistic case studies of hurricane-induced flooding in the Texas coastal region, one based on Hurricane Harvey and the other on Tropical Storm Imelda. We assess the impact of the chosen surrogate PF model, the uncertainty perspective, and the mitigation budget on the optimal flood mitigation. Our comparisons of surrogate PF models extend those presented in our preliminary assessments of fixed-topology grids in a winter storm resilience application \citep{Austgen2022a}, and our comparison of the adopted uncertainty perspective is similar to that presented in \citet{Shukla2023} for long-term flood mitigation. Our results suggest the adopted uncertainty perspective and given resource budget are impactful whereas choosing an LPAC model instead of a DC model for the recourse problem has limited effect. Because incorporating LPAC incurs a relatively large computational cost, we deduce that the DC approximation is preferred for this time-sensitive resilience application.
    \item To validate our models, we assess the performance of the mitigation solutions they prescribe in an ACPF model. Our results suggest that our models, despite their inaccuracies, are apparently capable of prescribing effective flood mitigation.
\end{enumerate}

The remainder of the paper is structured as follows. Section~\ref{section:review} provides a review of literature pertaining to resilience, flood modeling, and PF modeling. In Section~\ref{section:modeling}, we introduce our adaptations of the DC and LPAC power flow approximation models and their incorporation in overarching two-stage stochastic programming and robust optimization models designed to inform flood mitigation decisions for electric substations. In Section~\ref{section:results}, we compare the models using results obtained from applying the models to the flooding case studies. Finally, we present our conclusions in Section~\ref{section:conclusions}.
\section{Literature Review} \label{section:review}

In this paper, we introduce a short-term approach to combating substation flooding induced by an imminent but still-uncertain hurricane, a high-impact, low-frequency (HILF) event. The scenarios available for contingency planning are highly dependent on the near-term weather forecasts and thus strongly correlated. In this section, we discuss approaches to similar problems in the literature and their connection to our application.

\subsection{Reliability, Resilience, and Risk}

Mid- and long-term power grid contingency planning has historically revolved around reliability indices like System Average Interruption Frequency Index (SAIFI), System Average Interruption Duration Index (SAIDI), and others from the IEEE Guide for Electric Power Distribution Reliability Indices standards document \citep{IEEEReliabilityIndices2022}. These look-behind, distribution system-oriented metrics are inadequate for anticipating the health of transmission systems following geographically broad extreme weather events like hurricanes. While the security-constrained optimal power flow (SCOPF) model has been useful for proactive planning, it has typically been used to mitigate contingency sets like $N-1$ \citep{Zhang2012,Dvorkin2018} or $N-k$ \citep{Huang2022}, which are reliability-focused and not generally representative of the severe, correlated states that characterize the uncertainty of of an imminent threat. Additionally, the purpose of SCOPF is to ensure an acceptable or ``secure'' outcome. Extreme weather events like hurricanes often engender consequences too severe to ensure ``secure'' grid operation, especially when mitigation resources are limited.

To overcome these limitations, we adopt an approach based on risk and resilience principles. From a 2012 report by the National Research Council \citep{NRC2012}, ``resilience is the ability to prepare and plan for, absorb, recover from, and more successfully adapt to adverse events.'' More nuanced aspects of risk and resilience are debated. For example, risk and resilience as defined in \citet{Linkov2019} are respectively threat-dependent and threat-agnostic. However, in \citet{Logan2022}, resilience is defined as being system- and context-specific and inherently integrated with risk. We adopt the latter perspective and develop our model according to the conceptual framework in \citet{Watson2014} that proposes resilience metrics be formed with consideration of three factors: the threat, the likelihood, and the consequences. In our models, we capture the threat and likelihoods by considering a sample of representative hurricane flooding scenarios. Consequences are captured by an optimal power flow (OPF) model that is a function of both the flooding realization and preparedness decisions. To capture risk, we employ as risk measures the expectation operator associated with stochastic programming (SP) \citep{Birge2011} and the maximum operator associated with robust optimization (RO) \citep{BenTal2009}.

\subsection{Flood Mitigation and Power Flow Modeling}

Though research on substation flood mitigation is limited, research on general flood mitigation is extensive. For example, \citet{Eijgenraam2014,Eijgenraam2017} propose an optimization model for planning dike ring heightening in the Netherlands. \citet{Klerk2021} propose a heuristic method for solving a similar problem that considers dike ring heterogeneity. For mitigating specific flooding events (\textit{e.g.}, a dam break), \citet{Tasseff2019} introduce an approach based on a partial differential equation model of flooding dynamics that incorporates wall-building and revegetation as mitigation decisions. Though all of these models are more sophisticated than ours in some ways, none of them consider the dynamics of the at-risk infrastructure.

In power grid resilience research, there are numerous examples of surrogate PF models being incorporated in optimization models to capture such dynamics. \citet{Movahednia2022a} propose a network-agnostic substation importance index for the highly related problem of scheduling \tigerdam~deployments prior to an imminent hurricane. Examples of network-based surrogate power flow models include
the network flow relaxation \citep{Tan2018,Mohagheghi2015},
the DC power transfer distribution factor (PTDF) approximation \citep{Garcia2022},
the DC B-theta approximation \citep{Arab2015,Movahednia2022b,Shukla2023,Souto2022,Yang2023,Pierre2018,Coffrin2011,SahraeiArdakani2017,Quarm2022},
the LPAC approximation \citep{Coffrin2015},
and the second-order cone programming (SOCP) relaxation \citep{Garifi2022}.

Naturally, the complexity of a power system resilience problem depends on many factors including (i) the embedded power flow model, (ii) properties of the instance such as grid size, grid topology, uncertainty characterization, etc., (iii) the scope of the resilience decisions (\textit{i.e.}, mitigation and response, mitigation and restoration, restoration only, etc.), and (iv) the power system controls considered. Regarding (iv), a significant driver of complexity is the inclusion of discrete, potentially topology-affecting controls (\textit{e.g.}, branch switching) in the model. To attain a computationally tractable instance, one must consider trade-offs among these factors. For example, \citet{Yang2023} incorporate discrete bus, branch, and generator de-energization decisions in a relatively simple multi-time period DC power flow model. The overarching decision problem is a two-stage stochastic programming model, which is applied to a case study featuring the 73-bus RTS-GMLC grid instance and 500 wildfire contingency scenarios. \citet{Wang2013} also use a multi-time period DC power flow model for their $N-k$ contingency-constrained unit commitment problem. In that model, discrete variables are used to capture exogenous generator commitment decisions and also to endogenously indicate the worst $N-k$ contingency. The problem is formulated as a two-stage robust optimization problem, which is evaluated on a modified IEEE 118-bus test system for $k = 1, 2, 3, 4$. In contrast, \citet{Garifi2022} leverage the relatively complex SOCP power flow model to inform binary bus, branch, and generator hardening decisions as well as contingency response decisions, including binary generator and branch switching decisions, over multiple time periods. The overarching problem is formulated as a single-stage deterministic optimization model, which is applied to the 73-bus RTS-GMLC grid instance. \citet{Coffrin2019} incorporate feasibility-ensuring discrete decisions for switching buses, bus shunts, and generators into the AC power flow model. The model is used to effectively evaluate the ``maximal load delivery'' response to a given contingency for grid instances with up to roughly 2000 buses. Similarly, \cite{Rhodes2021} incorporate discrete component re-energization decisions into the AC power flow model, as well as its SOCP relaxation and DC approximation. The models are used to find an optimal restoration plan over multiple time periods. The AC and SOCP formulations are applied to a 5-bus instance, and the DC formulation is applied to the IEEE 118-bus and 123-bus radial distribution test systems.

In this paper, we extend our past research on the deployment of flood mitigation resources prior to an imminent but still uncertain hurricane threat to minimize the load shed under the resulting contingency \citep{Austgen2023}. To capture the operational response, we consider the two-stage SP model with recourse based on DC power flow presented in that paper and additionally two-stage RO formulations and variants of the relatively more complex LPAC approximation in the recourse model. In \cite{Austgen2022a}, we conducted a preliminary study on fixed-topology forms of these power flow models in the context of winter storm planning. In this paper, we extend the models to account for mitigation- and contingency-dependent topology as in \cite{Pierre2018,Movahednia2022b,Garifi2022}. Due to the relatively complex scope of our resilience decisions (\textit{i.e.}, mitigation decisions under uncertainty) and the relatively large 663-bus grid instance we study, the only discrete exogenous decisions we incorporate in our model are the first-stage mitigation decisions. However, unlike in our past models, we now also incorporate a feasibility-guaranteeing endogenous indicator variable in the recourse problem that signals whether power system controls not captured by our model must be exercised to attain a feasible power flow.

\subsection{Power Flow Model Comparisons} \label{subsection:pf_model_comparisons}

For decision-making applications, choosing a surrogate PF model from the many that exist is a difficult task -- incorporating a more fidelitous PF model may sacrifice computational tractability while incorporating a less fidelitous PF model may compromise the overarching optimization model's ability to prescribe good decisions. This trade-off has been studied for a variety of power systems applications. For example, \citet{Overbye2004} concluded that the DC model yields locational marginal prices (LMPs) comparable to those yielded by the AC model, though it should be noted that LMPs are not decisions in those models but rather values computed from the obtained power flow solution. Other examples include the DC and LPAC models for system operation \citep{Coffrin2014} and restoration \citep{Coffrin2015}, DC and AC models for system reliability \citep{Kile2014}, and DC and SOCP models for transmission system hardening and restoration \citep{Garifi2022}. In these studies, the authors concluded the solutions of the DC model are typically worse than those of the higher-fidelity alternative in the context of an AC model feasibility study. For comparisons of decisions that are intricately linked to the PF model, \textit{e.g.}, the operational decisions in \citet{Coffrin2014}, such a conclusion may be expected. For comparisons of decisions more distant from the PF model, \textit{e.g.}, the mitigation investments in \citet{Garifi2022}, the effect of the chosen PF model is perhaps more difficult to anticipate.

In our past work, we compared the DC and LPAC models in a similar winter storm mitigation context and observed that that the mitigation decisions are largely unimpacted by the choice of these two power flow models; however, grid topology was assumed to be static in that study \citep{Austgen2022a}. In this paper, we similarly compare the DC and LPAC models in the context of flood mitigation planning, an application in which grid topology varies based on the mitigation decisions and flooding realizations.

\section{Modeling} \label{section:modeling}

\newcommand{\pathmodeling}{modeling}

\begin{figure}
    \FIGURE
    {\includegraphics[width=3in]{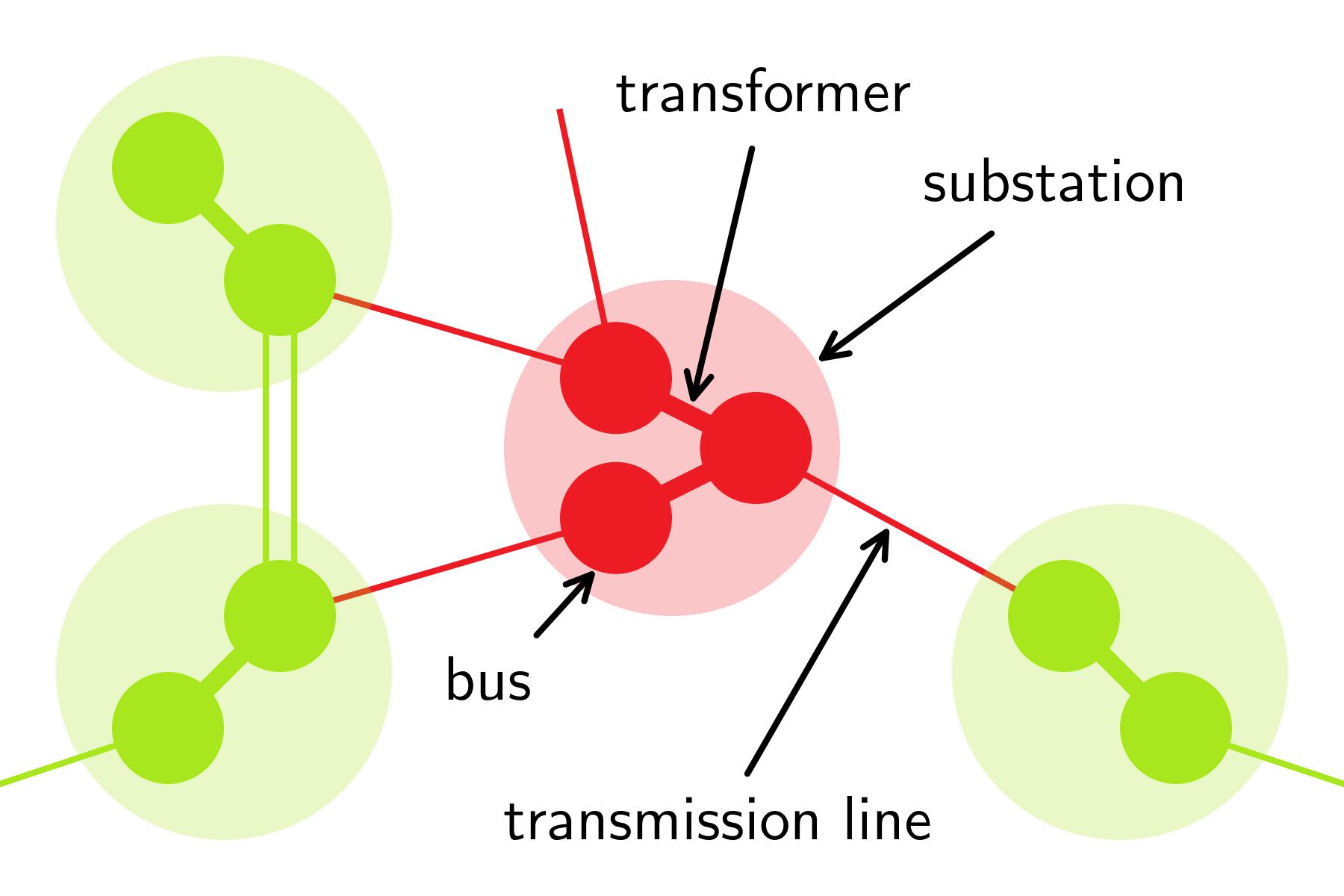}}
    {Impact of Substation Flooding on Buses and Branches \label{fig:grid_nomenclature}}
    {Here, the light red substation is flooded, and the red components are consequently inoperable. Light green substations are not flooded and the green components are operable.}
\end{figure}

For the sake of modeling, the power grid may be viewed as a graph with buses as nodes and branches (\textit{i.e.}, transmission lines and transformers) as edges. When a substation floods, we model all its buses and transformers, incident transmission lines, and associated generators as inoperable, and all its associated loads as unsatisfiable. This is illustrated in Figure~\ref{fig:grid_nomenclature}. We suppose at each substation that a discrete set of resilience levels are available for implementation by stacking temporary flood barriers like \tigerdams, and we use a PF model to assess the load shed that results in each scenario.

\subsection{Notation} \label{subsection:notation}

We now introduce the sets, parameters, and decision variables used in our models.

% some things to know about enumitem
% 1. spacing options are as described in the guide:
%    https://ctan.math.illinois.edu/macros/latex/contrib/enumitem/enumitem.pdf
% 2. there is some weirdness when using square brackets in math mode in an enumitem label;
%    this is remedied by putting the bracketed expression inside curly braces
% 3. the width of the label can be calculated via `calc` package:
%    https://tex.stackexchange.com/questions/18576/get-width-of-a-given-text-as-length

\noindent\textbf{Sets}
\begin{itemize}[align=left,leftmargin=10em,labelwidth=9em,labelsep=1em,itemsep=0em]
\item[$K$] set of substations
\item[$R = \{0, 1, \ldots, \hat{r}\}$] set of resilience levels
\item[$\Omega$] set of scenarios
\item[$N$, $N_k$] set of buses, buses at substation $k$
\item[$E$] set of edges (in LPAC, $E = E^T \cup E^F$, the union of `to' and `from' edges)
\item[$G$, $G_n$] set of generators, generators at bus $n$
\item[$D$, $D_n$] set of loads, loads at bus $n$
\item[$L$, $L_{nm}$] set of branches, branches defined on edge $(n,m)$
\item[$L_n^-$, $L_n^+$, $L_n$] set of branches toward node $n$, away from node $n$, incident to node $n$
\item[$\Theta_{\text{cos}} \subseteq {[-\pi, \pi]}$] set of tangent line intersection points in the polyhedral relaxation of cosine
\item[$\Theta_\text{disc} \subseteq {[0, 2\pi]}$] set of tangent line intersection points in the polyhedral relaxation of a disc
\end{itemize}

\noindent\textbf{Parameters}
\begin{itemize}[align=left,leftmargin=10em,labelwidth=9em,labelsep=1em,itemsep=0em]
\item[$\prob(\omega)$] probability of scenario $\omega$ occurring
\item[$\lambda^\text{shed}$, $\lambda^\text{over}$] unitless objective weights for load shed and overgeneration
\item[$c_{kr}$] marginal resource cost of reinforcing substation $k$ to resilience level $r$ from level $r-1$
\item[$f$] resource budget
\item[$\xi_{kr}$] flooding uncertainty; 1 if substation $k$ is flooded to level $r$; 0 otherwise
\item[$b_l$, $g_l$] susceptance and conductance of branch $l$
\item[$\underline{p}_g^\text{gen}, \underline{q}_g^\text{gen}$] lower bounds for active and reactive power generation at generator $g$
\item[$\overline{p}_g^\text{gen}, \overline{q}_g^\text{gen}$] upper bounds for active and reactive power generation at generator $g$
\item[$p_d^\text{load}$, $q_d^\text{load}$] active and reactive power demand of load $d$
\item[$\overline{s}_l^\text{flow}$] upper bound of apparent power flow across branch $l$
\item[$\underline{v}_n, v_n, \overline{v}_n$] minimum, target, and maximum voltage magnitudes of bus $n$
\item[$n_\text{ref}$] reference bus
\item[$\overline{\theta}_\Delta$] maximum voltage phase angle difference of adjacent buses
\item[$\overline{\theta}$] maximum absolute voltage phase angle for any bus in the system
\item[$M$] arbitrarily large positive constant (for big-$M$ method)
\end{itemize}

\noindent\textbf{Decision Variables}
\begin{itemize}[align=left,leftmargin=10em,labelwidth=9em,labelsep=1em,itemsep=0em]
\item[$x_{kr} \in \{0,1\}$] flood mitigation indicator variable; 1 if substation $k$ is resilient to level $r$, 0 otherwise
\item[$\chi \in \{0,1\}$] infeasibility indicator variable; 1 if no power flow solution exists, 0 otherwise
\item[$\alpha_n \in \{0,1\}$] bus status indicator variable; 1 if bus $n$ is operational, 0 otherwise
\item[$\beta_{nm} \in \{0,1\}$] branch status indicator variable; 1 if branches on edge $(n,m)$ are operational, 0 otherwise
\item[$\hat{p}_g \in \mathbb{R}$, $\hat{q}_g \in \mathbb{R}$] active and reactive power generation at generator $g$
\item[$\check{p}_g \in \mathbb{R}_+$] active power overgeneration at generator $g$
\item[$\tilde{p}_l \in \mathbb{R}$, $\tilde{q}_l \in \mathbb{R}$] active and reactive power flow across branch $l$
% \delta_d was formerly z_d, but I am now using z to denote objective value in results chapter
\item[$\delta_d \in {[0,1]}$] proportion of load $d$ satisfied
\item[$\theta_n \in {[-\overline{\theta},\overline{\theta}]}$] voltage phase angle of bus $n$
\item[$\widehat{\sin}_{nm} \in \mathbb{R}$] approximation of $\sin(\theta_n - \theta_m)$
\item[$\widehat{\cos}_{nm} \in {[\cos(\overline{\theta}_\Delta), 1]}$] approximation of $\cos(\theta_n - \theta_m)$
\item[$\phi_n \in \mathbb{R}$] deviation from target voltage magnitude at bus $n$
\end{itemize}

We aim to denote similar types of parameters and variables similarly. For example, all power flow variables have a tilde (\textit{e.g.}, $\tilde{p}_l$ and $\tilde{q}_l$), and all power generation variables have a hat (\textit{e.g.}, $\hat{p}_g$ and $\hat{q}_g$). Underlining and overlining signify lower and upper bounds, respectively (\textit{e.g.}, $\underline{p}_g^\text{gen}$ and $\overline{p}_g^\text{gen}$). Additionally, when a parameter or variable name appears in bold, it denotes the vector comprising all the indexed elements (\textit{e.g.}, $\bs{x} = [x_{kr}, \forall k \in K, \forall r \in R]$). When necessary, we apply a superscript $\omega$ (\textit{e.g.}, $\bs{\xi}^\omega$) to indicate a quantity associated specifically with scenario $\omega$. All power grid parameters and variables are assumed to be in the per-unit system.
\subsection{Resilience Modeling}

The discrete set of implementable resilience levels is denoted by $R$, and the decision to reinforce substation $k$ to a specific level $r$ is captured by the binary decision variable $x_{kr}$. Substation $k$ is resilient to level $r$ flooding if $x_{kr} = 1$ and is otherwise susceptible. The mitigation model comprises three sets of constraints involving these variables:
\begin{subequations}
\begin{gather}
    x_{k,r+1} \le x_{k,r},\quad \forall k \in K, \forall r \in R \setminus \{\hat{r}\},
    \label{eq:con_incremental} \\
    x_{k\hat{r}} = 0,\quad \forall k \in K,
    \label{eq:con_inexorable} \\
    \sum_{k \in K} \sum_{r \in R} c_{kr} x_{kr} \le f.
    \label{eq:con_resource_hi}
\end{gather}
\end{subequations}
Constraints \eqref{eq:con_incremental} capture the cumulative nature of mitigation. We suppose the mitigation is limited either physically or practically and introduce $\hat{r} = \max\{R\}$ as the unattainable level of resilience. Flooding at or above the unattainable level is rendered inexorable by constraints \eqref{eq:con_inexorable}. Finally, we suppose the mitigation resources are limited to a budget of $f$.
Supposing a marginal cost $c_{kr}$ associated with each decision $x_{kr}$, this is captured by the binary knapsack constraint \eqref{eq:con_resource_hi}. For brevity, we hereafter refer to the constraints of the mitigation decision making problem as
\begin{equation*}
    \mathcal{X} = \{\boldsymbol{x} \in \{0,1\}^{|K \times R|} :
                    \eqref{eq:con_incremental},
                    \eqref{eq:con_inexorable},
                    \eqref{eq:con_resource_hi}\}. \label{eq:resilience_X}
\end{equation*}
\subsection{Power Flow Modeling} \label{subsection:power_flow_modeling}

Our recourse problems are based on two linear approximations of ACPF: the classical DC approximation \citep{Molzahn2019} and the more contemporary and fidelitous linear programming AC (LPAC) approximation \citep{Coffrin2014}.
In this section, we present adaptations of the models that account for incapacitated substations and incorporate various levels of detail.
In each model, we assume that a bus is operational if and only if its substation is operational. Similarly, we assume that a branch is operational if and only if both of its bus endpoints are operational.
Because these conditions are subject to change as a result of unpredictable flooding and prior mitigation decision making, each recourse problem is modeled as a function of a specific set of flood mitigation decisions and a specific flooding realization.

\subsubsection{AC Power Flow.}
Both the DC and LPAC approximations are based on the four core ACPF equations:
\begin{subequations}
\begin{alignat}{3}
& \sum_{g \in G_n} \hat{p}_g - \sum_{d \in D_n} p_n^\text{load} = \sum_{l \in L_n} \tilde{p}_l,
&& \quad \forall n \in N, \label{eq:acpf_p_kcl} \\
& \sum_{g \in G_n} \hat{q}_g - \sum_{d \in D_n} q_n^\text{load} = \sum_{l \in L_n} \tilde{q}_l,
&& \quad \forall n \in N, \label{eq:acpf_q_kcl} \\
& \tilde{p}_l = v_n^2 g_l - v_n v_m g_l \cos(\theta_n - \theta_m) - v_n v_m b_l \sin(\theta_n - \theta_m),
&& \quad \forall (n,m) \in E, \forall l \in L_{nm}, \label{eq:acpf_p_ohms_law} \\
& \tilde{q}_l = -v_n^2 b_l + v_n v_m b_l \cos(\theta_n - \theta_m) - v_n v_m g_l \sin(\theta_n - \theta_m),
&& \quad \forall (n,m) \in E, \forall l \in L_{nm} \label{eq:acpf_q_ohms_law}.
\end{alignat}
\end{subequations}
These equations are based on the simple series admittance model of a branch \citep{Molzahn2019}. Here, complex power is represented in rectangular form (\textit{i.e.}, $p + jq$) and complex voltage in polar form (\textit{i.e.}, $ve^{j\theta}$). Linear equations \eqref{eq:acpf_p_kcl} and \eqref{eq:acpf_q_kcl} capture Kirchhoff's Current Law (KCL), and nonconvex equations \eqref{eq:acpf_p_ohms_law} and \eqref{eq:acpf_q_ohms_law} capture Ohm's Law. Many approximation, relaxation, and restriction PF models are obtained by manipulating the nonconvex equations to make them more computationally tractable.

\subsubsection{Adapted DC Model.}
% intuitive explanation of how DCOPF is derived:
% https://invenia.github.io/blog/2021/06/18/opf-intro/
The classical DC approximation is based on three assumptions: (1) branch conductance is negligible relative to susceptance and may be neglected, (2) bus voltage magnitudes may be fixed to one per unit, and (3) the difference between voltage phase angles of adjacent buses is small such that the sine of that difference may be accurately modeled by the difference itself and the cosine may be fixed to one.
These assumptions lead to the reactive power flows being zero and the active power flows obeying a linear relationship with the bus voltage phase angles.
This model may be viewed as an extension of a capacitated network flow problem involving multiple sources and sinks. Importantly, the solution space is additionally confined by the complicating Ohm's Law equality constraints.
Our adaptation of the DC approximation is
\begin{subequations} \label{eq:L_dc}
\begin{alignat}{2}
\mathllap{\mathcal{L}_\text{DC}(\bs{x}, \bs{\xi}) =}
% https://tex.stackexchange.com/a/324252/290597
% jank... note this is an issue in TeX 2020, but not in TeX 2022
\setlength{\dimen0}{\widthof{$\text{min}~~$}}
\hspace{\dimen0}
\setlength{\dimen0}{0pt-\widthof{$\mathcal{L}_\text{DC}(\bs{x}, \bs{\xi}) =$}}
\hspace{\dimen0}
\hspace{2em}
&&& \notag \\
\text{min}~~& \lambda^\text{shed} \sum_{d \in D} p_d^\text{load} (1 - \delta_d) + \lambda^\text{over} \sum_{g \in G} \check{p}_g
&& \label{eq:dc_min_load_shed} \\
% linking constraints
\text{s.t.}~~& \alpha_n = \prod_{r \in R} \left(1 - \xi_{kr} \left(1 - x_{kr}\right)\right),
&& \quad \forall k \in K, \forall n \in N_k, \label{eq:def_alpha} \\
& \beta_{nm} = \alpha_n \alpha_m,
&& \quad \forall (n,m) \in E, \label{eq:def_beta} \\
% load satisfaction prevention constraint
& \delta_d \le 1 - \chi,
&& \quad \forall d \in D, \label{eq:infeasibility_load_shed} \\
% flow conservation
& \sum_{g \in G_n} \left(\hat{p}_g - \check{p}_g\right) + \sum_{l \in L_n^-} \tilde{p}_l
&& \notag \\*
&
\quad = \sum_{d \in D_n} p_d^\text{load} \delta_d + \sum_{l \in L_n^+} \tilde{p}_l,
&& \quad \forall n \in N, \label{eq:dc_kcl} \\
% Ohm's Law
& M \left(1 - \beta_{nm}\right) \ge \left| -\tilde{p}_l - b_l \widehat{\sin}_{nm} \right|,
&& \quad \forall (n,m) \in E, \forall l \in L_{nm}, \label{eq:dc_ohms_law} \\
% sin approximations
& \widehat{\text{sin}}_{nm} \in \mathcal{Y}_\text{sin}(\theta_n, \theta_m, \beta_{nm}),
&& \quad \forall (n,m) \in E, \label{eq:dc_sin} \\
% branch limits
& -\overline{s}_l^\text{flow} \beta_{nm} \le \tilde{p}_l \le \overline{s}_l^\text{flow} \beta_{nm},
&& \quad \forall (n, m) \in E, \forall l \in L_{nm}, \label{eq:dc_flow_limits} \\
% generator limits
& \underline{p}_g^\text{gen} (\alpha_n - \chi) \le \hat{p}_g \le \overline{p}_g^\text{gen} \alpha_n,
&& \quad \forall n \in N, \forall g \in G_n, \label{eq:dc_generation_limits} \\
% overgeneration constraint
& \check{p}_g \le \hat{p}_g,
&& \quad \forall g \in G, \label{eq:dc_overgeneration_limit} \\
% reference bus voltage fixing
& \theta_{n_\text{ref}} = 0.
&& \label{eq:dc_ref_voltage_phase_angle}
\end{alignat}
\end{subequations}
Here, variables are also constrained as specified in Section~\ref{subsection:notation}. The objective \eqref{eq:dc_min_load_shed} is to minimize the weighted combination of active power load shed and overgeneration. Constraints \eqref{eq:def_alpha} and \eqref{eq:def_beta} relate the operational statuses of buses and branches to those of the substations. The operational status of each component is determined exactly by $\boldsymbol{x}$ and $\boldsymbol{\xi}$, the arguments to $\mathcal{L_\text{DC}}$. That is, the domain of each $\alpha_n$ and $\beta_{nm}$ may be relaxed to $[0, 1]$ since their integrality is guaranteed by constraints \eqref{eq:def_alpha} and \eqref{eq:def_beta} and the integrality of $\bs{x}$.
The logical formulation here is nonlinear, but as detailed in \citet{Asghari2022}, constraints \eqref{eq:def_alpha} may be linearly reformulated as
\begin{alignat*}{3}
    \alpha_n &\ge \displaystyle \sum_{r \in R} \left(1 - \xi_{kr} \left(1 - x_{kr}\right)\right) - |R| + 1,
    && \quad \forall k \in K, \forall n \in N_k, \\
    \alpha_n &\le 1 - \xi_{kr} \left(1 - x_{kr}\right),
    && \quad \forall k \in K, \forall n \in N_k, \forall r \in R,
\end{alignat*}
and constraints \eqref{eq:def_beta} as
\begin{alignat*}{3}
    \beta_{nm} &\ge \alpha_n + \alpha_m - 1, && \quad \forall (n,m) \in E, \\
    \beta_{nm} &\le \alpha_n,                && \quad \forall (n,m) \in E, \\
    \beta_{nm} &\le \alpha_m,                && \quad \forall (n,m) \in E.
\end{alignat*}
Constraints \eqref{eq:infeasibility_load_shed} ensure that no load is satisfied when $\chi = 1$. Constraints \eqref{eq:dc_kcl} impose KCL, the PF equivalent of flow balance. In the standard DC approximation, Ohm's Law is represented as the equality constraint ${-p_l - b_l \widehat{\sin} = 0}$. To ensure out-of-service branches are treated as open circuits, we embed the big-$M$ technique from \citet{Coffrin2011} in constraints \eqref{eq:dc_ohms_law} to enforce the equality only for operational branches. Details on how the big-$M$ constants may be calibrated are presented in Appendix B in the online supplement. The approximation of $\sin(\theta_n - \theta_m)$ is captured by constraints \eqref{eq:dc_sin}. In our implementation,
\begin{align}
    \mathcal{Y}_\text{sin}\left(\theta_n, \theta_m, \beta_{nm}\right)
    = \Big\{
        \widehat{\text{sin}} :
        &\!~\widehat{\text{sin}} = \theta_n - \theta_m, 
        \big| \widehat{\sin} \big| \le 2 (1 - \beta_{nm}) \overline{\theta} + \beta_{nm} \overline{\theta}_\Delta
    \Big\}. \label{eq:variant_sine_linear}
\end{align}
That is, we model $\sin(\theta_n - \theta_m) \approx \theta_n - \theta_m$ with $\widehat{\sin} \in [-\overline{\theta}_\Delta, \overline{\theta}_\Delta]$ if $\beta_{nm} = 1$ and $\widehat{\sin} \in [-2 \overline{\theta}, 2 \overline{\theta}]$ if $\beta_{nm} = 0$ so that both $\theta_n$ and $\theta_m$ are free to assume any value in their ordinary $[-\overline{\theta}, \overline{\theta}]$ range. Constraints \eqref{eq:dc_flow_limits} impose conditional lower and upper bounds on power flows. Power generation lower and upper bounds are imposed by constraints \eqref{eq:dc_generation_limits}. The bounds depend on the operational status of the corresponding bus and on the infeasibility indicator $\chi$. In Section~\ref{subsection:relatively_complete_recourse}, we discuss the role of this variable in guaranteeing relatively complete recourse. Constraints \eqref{eq:dc_overgeneration_limit} ensure no more power is discarded than is generated at each generator, and constraint \eqref{eq:dc_ref_voltage_phase_angle} ensures the voltage phase angle of the reference bus is zero.

\subsubsection{Adapted LPAC Model.}
Recall one of the assumptions underlying the DC approximation is that bus voltage magnitudes may be set to one per unit. This assumption is reasonable when the grid is in good health, \textit{e.g.}, when determining locational marginal prices \citep{Liu2009,Paul2017}. However, this assumption is less likely to hold when the grid is stressed by multiple damaged components as is often the case during and after extreme weather events.

In \citet{Coffrin2014}, hot-, warm-, and cold-start variants of LPAC are proposed for situations where all, some, or no information about bus voltage magnitudes are available. We leverage the warm-start LPAC approximation which incorporates deviations from target bus voltage magnitudes. This grants some flexibility in the feasible solutions, and stable grid operation is still ensured by bounding the deviations. The warm-start variant additionally incorporates reactive power and line losses. The approximation is based on four assumptions, some the same as for the DC approximation: (1) the difference of adjacent buses' voltage phase angles is small such that the sine of that difference may be accurately modeled by the difference itself, (2) the cosine of adjacent buses' voltage phase angles may be accurately modeled by a polyhedral relaxation of cosine, (3) the effects of bus voltage magnitude deviations on active power flow are negligible, and (4) target bus voltage magnitudes are available, and the effects of bus voltage magnitude deviations on reactive power flow may be modeled by a linear approximation of bus voltage around the target. The LPAC approximation is developed using variable substitutions, first-order Taylor polynomial approximations, and McCormick envelopes. Just as the DC approximation may be interpreted as a sort of network flow problem, the LPAC approximation may be viewed as a capacitated multi-source, multi-sink, and moreover multi-layer network flow problem in which intra- and inter-network interactions are dictated by the complicating Ohm's Law equality constraints. Our adaptation of this model is
\begin{subequations} \label{eq:L_lpac}
\begin{alignat}{2}
\mathllap{\mathcal{L}_\text{LPAC}(\bs{x}, \bs{\xi}) =}
% https://tex.stackexchange.com/a/324252/290597
% jank... note this is an issue in TeX 2020, but not in TeX 2022
\setlength{\dimen0}{\widthof{$\text{min}~~$}}
\hspace{\dimen0}
\setlength{\dimen0}{0pt-\widthof{$\mathcal{L}_\text{LPAC}(\bs{x}, \bs{\xi}) =$}}
\hspace{\dimen0}
\hspace{2em}
&&& \notag \\
 \text{min}~~& \lambda^\text{shed} \sum_{d \in D} p_d^\text{load} (1 - \delta_d) + \lambda^\text{over} \sum_{g \in G} \check{p}_g
&& \label{eq:lpac_min_load_shed} \\
% linking constraints
\text{s.t.}~~& \eqref{eq:def_alpha}, \eqref{eq:def_beta}, \eqref{eq:infeasibility_load_shed}, && \\
% flow conservation
& \sum_{g \in G_n} \left(\hat{p}_g - \check{p}_g\right) = \sum_{d \in D_n} p_d^\text{load} \delta_d + \sum_{l \in L_n} \tilde{p}_l,
&& \quad \forall n \in N, \label{eq:lpac_p_kcl} \\
& \sum_{g \in G_n} \hat{q}_g = \sum_{d \in D_n} q_d^\text{load} \delta_d + \sum_{l \in L_n} \tilde{q}_l,
&& \quad \forall n \in N, \label{eq:lpac_q_kcl} \\
% Ohm's Law
%
& M \left(1 - \beta_{nm}\right) \ge \Big| -\tilde{p}_l + v_n g_l \left(v_m - v_n\right) \chi + v_n^2 g_l && \notag \\*
& \phantom{M \left(1 - \beta_{n,m}\right) \ge \Big|}
- v_n v_m \left(g_l \widehat{\text{cos}}_{nm} + b_l \widehat{\text{sin}}_{nm}\right) \Big|,
&& \quad \forall (n,m) \in E, \forall l \in L_{nm}, \label{eq:lpac_p_ohms_law} \\
& M \left(1 - \beta_{nm}\right) \ge \Big| -\tilde{q}_l + v_n b_l (v_n - v_m) \chi - v_n^2 b_l && \notag \\*
& \phantom{M \left(1 - \beta_{n,m}\right) \ge \Big|}
- v_n v_m \left(g_l \widehat{\text{sin}}_{nm} - b_l \widehat{\text{cos}}_{nm}\right) && \notag \\*
& \phantom{M \left(1 - \beta_{n,m}\right) \ge \Big|}
- v_n b_l \left(\phi_n - \phi_m\right) - (v_n - v_m) b_l \phi_n \Big|,
&& \quad \forall (n,m) \in E, \forall l \in L_{nm}, \label{eq:lpac_q_ohms_law} \\
% cos and sin approximations
& \widehat{\text{sin}}_{nm} \in \mathcal{Y}_\text{sin}(\theta_n, \theta_m, \beta_{nm}),
&& \quad \forall (n,m) \in E, \label{eq:lpac_sin} \\
& \widehat{\text{cos}}_{nm} \in \mathcal{Y}_\text{cos}(\theta_n, \theta_m, \beta_{nm}),
&& \quad \forall (n,m) \in E, \label{eq:lpac_cos} \\
% branch limits
& (\tilde{p}_l, \tilde{q}_l) \in \mathcal{Y}_\text{disc}(\beta_{nm}),
&& \quad \forall (n, m) \in E, \forall l \in L_{nm}, \label{eq:lpac_branch_limits} \\
% generator limits
& \underline{p}_g^\text{gen} (\alpha_n - \chi) \le \hat{p}_g \le \overline{p}_g^\text{gen} \alpha_n,
&& \quad \forall n \in N, \forall g \in G_n, \label{eq:lpac_p_generation_limits} \\
& \underline{q}_g^\text{gen} \alpha_n \le \hat{q}_g \le \overline{q}_g^\text{gen} \alpha_n,
&& \quad \forall n \in N, \forall g \in G_n, \label{eq:lpac_q_generation_limits} \\
% overgeneration constraint
& \check{p}_g \le \hat{p}_g,
&& \quad \forall g \in G, \label{eq:lpac_p_overgeneration_limits} \\
% bus voltage magnitude limits
& \underline{v}_n \le v_n + \phi_n \le \overline{v}_n,
&& \quad \forall n \in N, \label{eq:lpac_bus_voltage_mag_limits} \\
% reference bus voltage fixing
& \theta_{n_\text{ref}} = 0. && \label{eq:lpac_ref_voltage_phase_angle}
\end{alignat}
\end{subequations}
Here, variables are also constrained as specified in Section~\ref{subsection:notation}.
The objective \eqref{eq:lpac_min_load_shed} and constraints \eqref{eq:def_alpha}, \eqref{eq:def_beta}, and \eqref{eq:infeasibility_load_shed} all serve the same purpose as they did in our adaptation of the DC approximation. Constraints \eqref{eq:lpac_p_kcl} and \eqref{eq:lpac_q_kcl} capture KCL by enforcing flow balance for both active and reactive power. Constraints \eqref{eq:lpac_p_kcl} in this model are similar but not exactly the same as constraints \eqref{eq:dc_kcl}. This is because each branch in the LPAC approximation is modeled by ``to'' and ``from'' flow variables, whereas each branch in the DC approximation is modeled by a single undirected flow variable. Constraints \eqref{eq:lpac_sin} and \eqref{eq:lpac_cos} model the sine and cosine of the difference of adjacent bus voltage phase angles, and constraints \eqref{eq:lpac_branch_limits} capture generally the thermal limit on apparent power flow. Constraints \eqref{eq:lpac_p_ohms_law} and \eqref{eq:lpac_q_ohms_law} enforce Ohm's Law for operational branches. Generator active and reactive power limits are enforced by \eqref{eq:lpac_p_generation_limits} and \eqref{eq:lpac_q_generation_limits}, and the active power overgeneration limits are enforced by \eqref{eq:lpac_p_overgeneration_limits}. Constraints \eqref{eq:lpac_bus_voltage_mag_limits} bound deviations from the target bus voltage magnitudes. Lastly, constraint \eqref{eq:lpac_ref_voltage_phase_angle} fixes the phase angle of the reference bus voltage.

For our adapted LPAC model, we develop three variants: a relatively coarse linear variant (LPAC-C), a relatively fine linear variant (LPAC-F), and a convex quadratic programming variant (QPAC). All three variants implement $\mathcal{Y}_\text{sin}$ as defined in \eqref{eq:variant_sine_linear}, but $\mathcal{Y}_\text{cos}$ and $\mathcal{Y}_\text{disc}$ are implemented differently. In \citet{Coffrin2014}, polygonal relaxations are proposed as a means of linearizing disc and cosine geometries. Of course, tighter polygons are preferred for maintaining an accurate model, but looser polygons are preferred for keeping the problem tractable. To study this trade-off, we formulate three variants of the LPAC model by incorporating different forms of $\mathcal{Y}_\text{disc}$ and $\mathcal{Y}_\text{cos}$.

\begin{figure}
    \FIGURE
    {
        \ifdefined\epsswitch
            \includegraphics{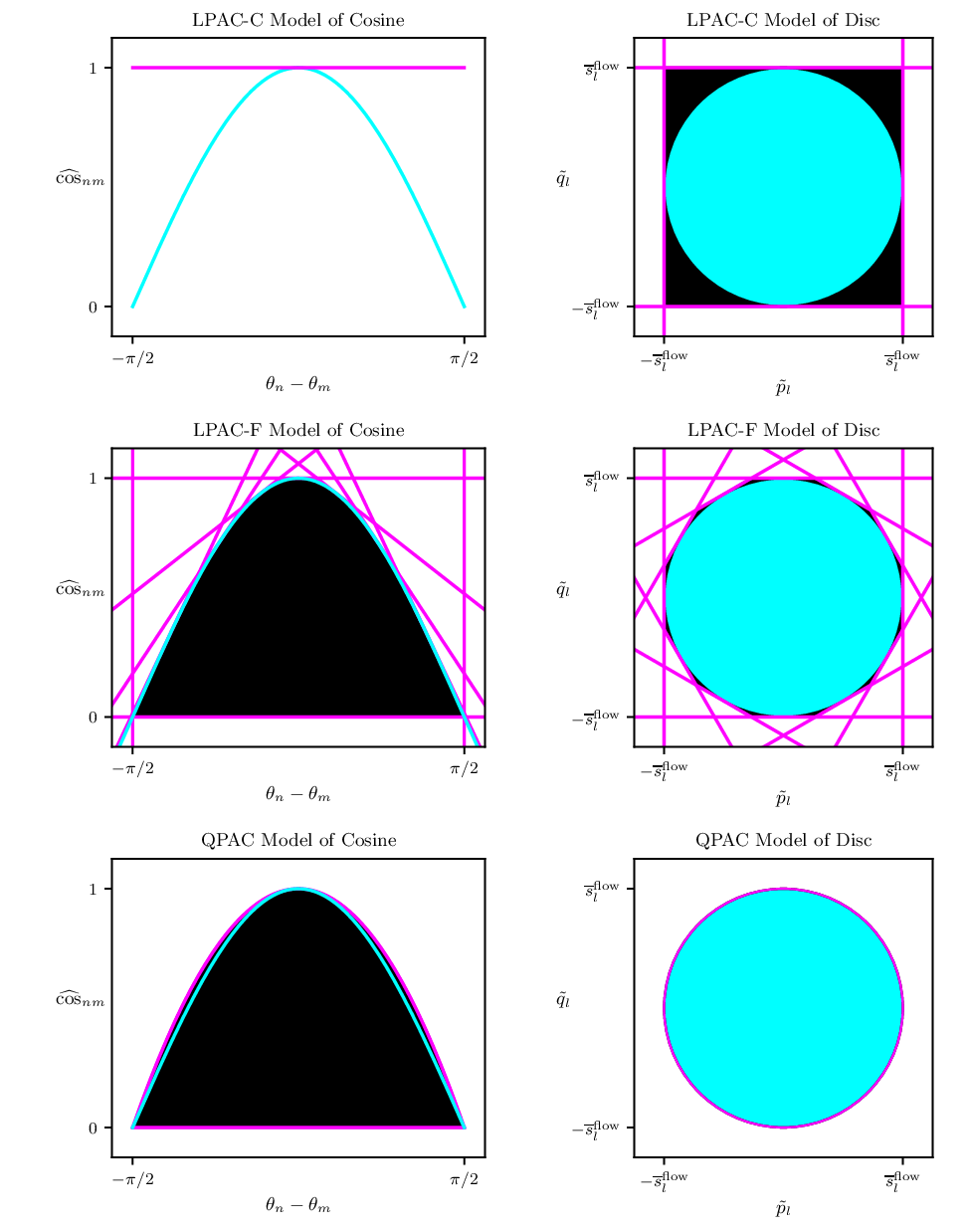}
        \else
            \includegraphics{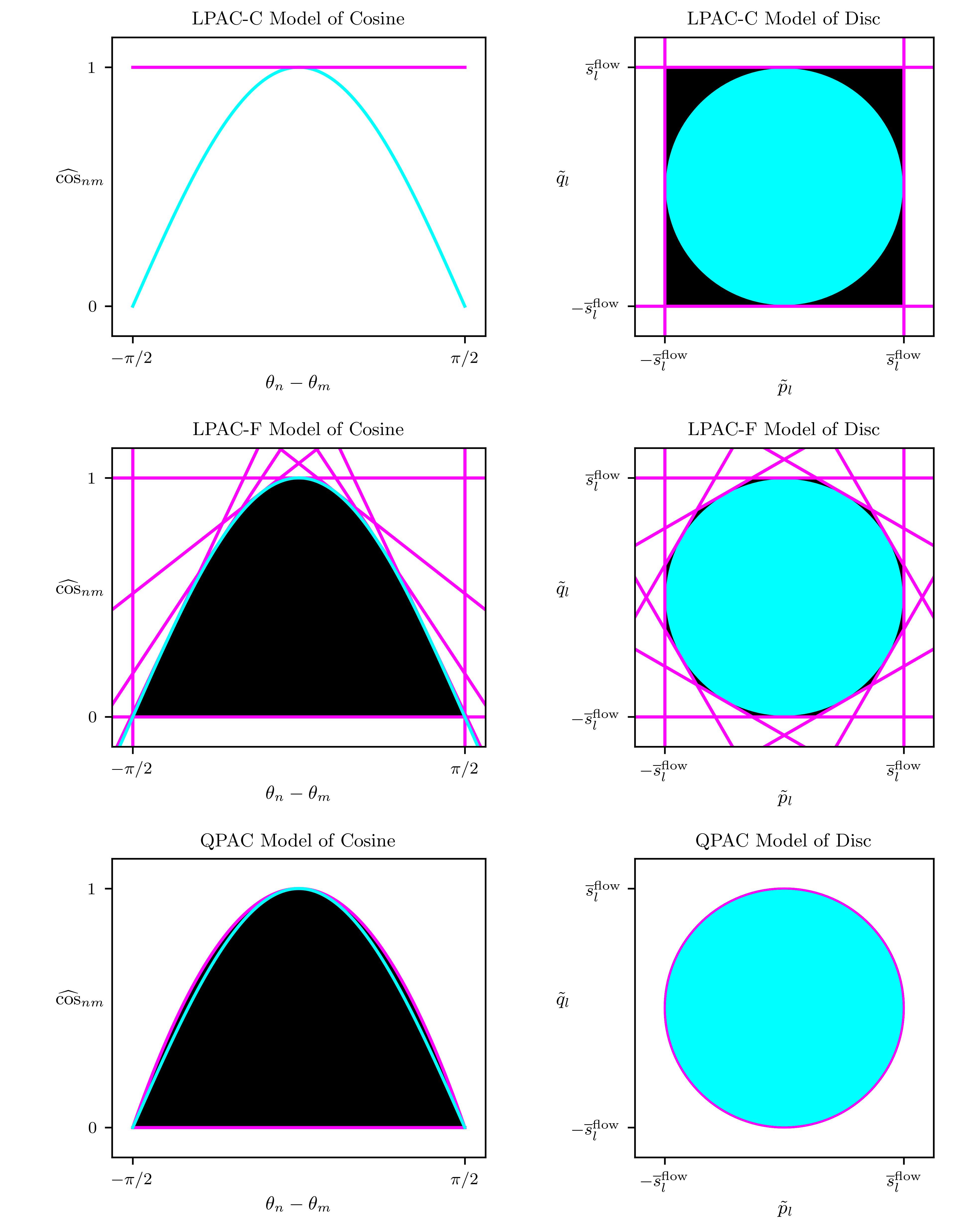}
        \fi
    }
    {Cosine and Disc Geometries for Our Three LPAC Variants for the Case of $\beta_{nm} = 1$ \label{fig:lpac_geometries}}
    {Cyan represents the exact sets, black the additional feasible solutions admitted by the relaxations, and magenta the boundaries of the linear and quadratic constraints that define the relaxed feasible regions.}
\end{figure}

In the LPAC-C variant, we adopt
\begin{gather*}
    \mathcal{Y}_\text{cos}\left(\theta_n, \theta_m, \beta_{nm}\right)
    = \left\{
        \widehat{\text{cos}} :
        \widehat{\text{cos}} = 1
    \right\}, \label{eq:lpac_l_cos} \\
     \mathcal{Y}_\text{disc}\left(\beta_{nm}\right)
    = \left\{
        \left(\tilde{p}, \tilde{q}\right) :
        \cos(\hat{\theta}) \tilde{p} + \sin(\hat{\theta}) \tilde{q} \le \overline{s}^\text{flow} \beta_{nm},~
        \forall \hat{\theta} \in \Theta_\text{disc}
    \right\} \label{eq:lpac_l_disc}
\end{gather*}
\noindent
with $\Theta_\text{disc} = \left\{\tfrac{t \pi}{2}, t = 1, \ldots, 4\right\}$ such that $\mathcal{Y}_\text{disc}$ is a square if $\beta_{nm} = 1$ and a set containing only $(0, 0)$ if $\beta_{nm} = 0$. This variant is based on one of the assumptions that underlies the DC approximation -- the difference of phase angles at adjacent buses is small such that the sine is approximately linear, and the cosine is approximately 1. For the case of $\beta_{nm} = 1$, these geometries are illustrated in the top row of subplots in Figure~\ref{fig:lpac_geometries}.

Our second, still linear, variant LPAC-F incorporates a polygonal relaxation of the disc and cosine. Define
\begin{equation*}
    \mathcal{B}_1(\theta; \hat{\theta}) = (\hat{\theta} - \theta) \sin(\hat{\theta}) + \cos(\hat{\theta})
\end{equation*}
as the line tangent to $\cos(\theta)$ at $\hat{\theta}$. In the LPAC-F variant, we let
\begin{equation*}
\begin{split}
    \mathcal{Y}_\text{cos}\left(\theta_n, \theta_m, \beta_{nm}\right)
    = \Big\{
        \widehat{\text{cos}} : \widehat{\text{cos}} \le 
        (1 - \beta_{nm}) (1 - \min\{\mathcal{B}_1(-2 \overline{\theta}; \hat{\theta}), \mathcal{B}_1(2 \overline{\theta}; \hat{\theta})\}) \quad \\
        + \mathcal{B}_1(\theta_n - \theta_m; \hat{\theta}),~
        \forall \hat{\theta} \in \Theta_\text{cos}
    \Big\} \label{eq:lpac_t_cos}
\end{split}
\end{equation*}
with ${\Theta_\text{cos} = \left\{0, \pm 0.354, \pm 0.735, \pm 1.211\right\}}$ such that $\mathcal{Y}_\text{cos}$ is a seven-edge polygonal relaxation of cosine when $\beta_{nm} = 1$. When $\beta_{nm} = 0$, the polygon is relaxed so that $\theta_n$ and $\theta_m$ are free to take values in ${[-\overline{\theta}, \overline{\theta}]}$ and $\widehat{\cos}$, an arbitrary decision variable in that case, is free to take values in ${[\cos(\overline{\theta}_\Delta), 1]}$. This specific $\Theta_\text{cos}$ is further motivated in Appendix A in the online supplement.
The LPAC-F variant also incorporates \eqref{eq:lpac_l_disc} but with ${\Theta_\text{disc} = \left\{\tfrac{t \pi}{6}, t = 1, \ldots, 12\right\}}$ such that $\mathcal{Y}_\text{disc}$ is a regular dodecahedron. For the case of $\beta_{nm} = 1$, these geometries are illustrated in the middle row of subplots in Figure~\ref{fig:lpac_geometries}.

Our third variant also retains the linear approximation of sine but employs the exact quadratic model of a disc and the tight convex quadratic relaxation of cosine from \citet{Molzahn2019}. In our work, we refer to this variant as the QPAC variant. Note that QPAC is sometimes used to refer to a specific quadratic programming ACPF approximation that was coincidentally patented by some of the same researchers that developed the LPAC approximation \citep{Coffrin2014,Coffrin2020}, but we do not mean it in that sense. Define
\begin{equation*}
    \mathcal{B}_2(\theta) = (1 - \cos(\overline{\theta}_\Delta)) \left(\theta / \overline{\theta}_\Delta\right)^2
\end{equation*}
as the convex quadratic curve that bounds $\cos(\theta)$ and intersects it at $-\overline{\theta}_\Delta, 0$, and $\overline{\theta}_\Delta$. Our QPAC variant implements
\begin{gather*}
    \mathcal{Y}_\text{cos}\left(\theta_n, \theta_m, \beta_{nm}\right)
    = \left\{
        \widehat{\text{cos}} :
        \widehat{\text{cos}} \le \mathcal{B}_2(\theta) + (1 - \beta_{nm}) (1 - \min\{\mathcal{B}_2(-2 \overline{\theta}), \mathcal{B}_2(2 \overline{\theta})\})
    \right\}, \label{eq:qpac_cos} \\
    \mathcal{Y}_\text{disc}\left(\beta_{nm}\right)
    = \left\{
        \left(\tilde{p}, \tilde{q}\right) :
        \tilde{p}^2 + \tilde{q}^2 \le (\overline{s}^\text{flow})^2 \beta_{nm}
    \right\}. \label{eq:qpac_disc}
\end{gather*}
Though this model sacrifices linearity, it allows the disc and cosine geometries to each be modeled using one quadratic inequality constraint (as opposed to multiple linear inequality constraints). In the case of $\beta_{nm} = 0$, these implementations behave the same as those from the LPAC-F variant. In the case of $\beta_{nm} = 1$, these geometries are illustrated in the bottom row of subplots in Figure~\ref{fig:lpac_geometries}.

\subsection{Two-Stage Models} \label{subsection:two_stage_models}

Having discussed the resilience and PF models, we now introduce our two-stage models. The first model seeks to minimize the expected load shed. We refer to it as the ``Stochastic Programming'' (SP) model, and it is formulated as
\begin{equation}
    \min_{\bs{x} \in \mathcal{X}} \sum_{\omega \in \Omega} \prob(\omega) \mathcal{L}(\bs{x}, \bs{\xi}^\omega)
    \tag{SP} \label{eq:SP}.
\end{equation}
Here, $\mathcal{X}$ is as defined in \eqref{eq:resilience_X} and $\mathcal{L}$ is either \eqref{eq:L_dc} or any variant of \eqref{eq:L_lpac}. In this model, the idea is to make first-stage resilience decisions with knowledge of how well the power grid will be able to perform in the aftermath of the hurricane flooding given those decisions. This model adopts the ``nature is fair'' approach -- each scenario $\omega$ is believed to occur with probability $\prob(\omega)$ and is weighted accordingly in the objective function.

Closely related to the SP model are a few solutions and bounds of theoretical and practical importance: the ``Expected Value`` (EV) solution, the ``Expected result of using the Expected Value solution'' (EEV) bound, and the ``Expected Wait-and-See'' (EWS) bound \citep{Birge2011}. An EV solution is obtained by solving a variant of the SP model in which the uncertainty is aggregated to its mean and is defined as ${\overline{\bs{x}} \in \argmin_{\bs{x} \in \mathcal{X}} \mathcal{L} \left( \bs{x}, \overline{\bs{\xi}} \right)}$.
For each scenario in the SP model, $\bs{\xi}^\omega = \mathcal{C}(\bs{\nu}^\omega)$ where ${\mathcal{C}: \mathbb{R}_{+}^{|K|} \to \{0, 1\}^{|K \times R|}}$ is a function that converts a real-valued vector $\bs{\nu}^\omega$ of flood levels to a vector $\bs{\xi}^\omega$ of flood indicators. Here, $\overline{\bs{\xi}} = \mathcal{C}(\overline{\bs{\nu}})$ where $\overline{\bs{\nu}} = \sum_{\omega \in \Omega} \prob(\omega) \bs{\nu}^\omega$. That is, our notion of the mean scenario is the set of flood indicators corresponding to the mean flood levels. The EEV bound is a function of the EV solution:
\begin{equation}
    \sum_{\omega \in \Omega} \prob(\omega) \mathcal{L}(\overline{\bs{x}}, \bs{\xi}^\omega)
    \tag{EEV} \label{eq:EEV}.
\end{equation}
The EEV bound captures how one would fare if one planned for the uncertainty's mean rather than its probability distribution. It is formulated as a restriction of the SP model that fixes the first-stage solution to an EV solution; thus, the EEV bound is an upper bound on the objective value of the SP model. The difference between the EEV bound and the objective value of the SP model is the ``Value of the Stochastic Solution'' (VSS).

The EWS bound is formulated similarly to the SP model but rather lets the first-stage decisions $\bs{x}$ be postponed until after the uncertainty is realized. It is formulated as
\begin{equation}
    \sum_{\omega \in \Omega} \prob(\omega) \min_{\bs{x} \in \mathcal{X}} \mathcal{L}(\bs{x}, \bs{\xi}^\omega).
    \tag{EWS} \label{eq:EWS}
\end{equation}
Consider that the SP model may be equivalently formulated by duplicating the $\bs{x}$ variables in each scenario and forcing the duplicate variables to be equal via non-anticipativity constraints. The EWS bound is formulated as the relaxation of the SP model that omits the nonanticipativity constraints. Thus, the EWS bound is a lower bound on the objective value of the SP model. The difference between the objective value of the SP model and the EWS bound is the ``Expected Value of Perfect Information'' (EVPI).

To contrast the SP model, we also propose a ``Robust Optimization'' (RO) model:
\begin{equation}
    \min_{\bs{x} \in \mathcal{X}} \max_{\omega \in \Omega} \mathcal{L}(\bs{x}, \bs{\xi}^\omega)
    \tag{RO} \label{eq:RO}.
\end{equation}
This model is much the same as the SP model, but it instead adopts the ``nature is adversarial'' approach. The probability of each scenario, if that information is even available, is ignored. Instead, the associated perspective is that the uncertainty realization will be that which does the most harm given the first-stage resilience decisions. This approach is often seen as overly if not impractically conservative. However, it is still valuable when scenario probabilities are not known or when the impact of one scenario dwarfs those of all other scenarios.

To accompany the RO model, we develop mathematical analogues to the EV solution, EEV bound, and EWS bound that we dub the ``Maximum Value'' (MV) solution, ``Maximum result of the Maximum Value solution'' (MMV) bound, and the ``Maximum Wait-and-See'' (MWS) bound. An MV solution is obtained by solving a variant of the RO model in which the uncertainty is aggregated to the substation-wise maximum and is defined as ${\widehat{\bs{x}} \in \argmin_{\bs{x} \in \mathcal{X}} \mathcal{L}(\bs{x}, \widehat{\bs{\xi}})}$. Here, $\widehat{\bs{\xi}} = \mathcal{C}(\widehat{\bs{\nu}})$ where $\widehat{\bs{\nu}} = \max_{\omega \in \Omega} \bs{\nu}^\omega$. That is, $\widehat{\bs{\nu}}$ is the scenario in which each substation experiences the worst flooding across any of the original scenarios, and we define the maximum value scenario as the corresponding set of flood indicators. Just as the RO model is a conservative approach to dealing with uncertainty, the MV solution is based on a conservative approach to constructing an aggregate scenario.

Naturally, the MMV bound then captures how one would fare in the RO model if one planned only for the maximum value scenario. Just like the EEV bound is formulated as a restriction of the SP model, the MMV bound is formulated as a restriction of the RO model and thus serves as an upper bound:
\begin{equation}
    \max_{\omega \in \Omega} \mathcal{L}(\widehat{\bs{x}}, \bs{\xi}^\omega)
    \tag{MMV} \label{eq:MMV}.
\end{equation}

Lastly, the MWS bound is formulated similarly to the RO model but rather lets the first-stage decisions $\bs{x}$ be postponed until after the uncertainty is realized. It is formulated as
\begin{equation}
    \max_{\omega \in \Omega} \min_{\bs{x} \in \mathcal{X}} \mathcal{L}(\bs{x}, \bs{\xi}^\omega).
    \tag{MWS} \label{eq:MWS}
\end{equation}
In the same way that the EWS bound is formulated as a relaxation of the SP model, the MWS bound is formulated as a relaxation of the RO model and thus serves as a lower bound.

Because the RO model is agnostic to the probability distribution of the scenarios (if one even exists), the mathematical analogues of the VSS and EVPI are not philosophically consistent with the perspective that nature is adversarial. As such, the MMV bound ought to be viewed simply as the performance that results from implementing $\widehat{\bs{x}}$. Similarly, the MWS bound ought to be viewed simply as a computationally inexpensive lower bound on the RO model's objective value.

\subsection{Relatively Complete Recourse} \label{subsection:relatively_complete_recourse}
When analyzing power grid contingencies using a PF model, it can occur that the model does not admit a feasible solution for certain contingencies. This is particularly troublesome for two-stage models like the SP and RO models because infeasibility in one or more scenarios causes the entire two-stage model to be infeasible. We design our adaptations of the DC and LPAC models with careful consideration of their role as recourse models in our two-stage formulations, and make certain modeling decisions to ensure their feasibility. Mainly, we incorporate a relaxed model of generator active power injections and also introduce an infeasibility indicator variable that grants relatively complete recourse \citep{Birge2011} by admitting a trivial ``shed all load'' solution when no non-trivial solution is feasible.

In power flow models, it is typical to model active  generation limits for a generator $g$ as
\begin{equation*}
    \underline{p}_g^\text{gen} \le \hat{p}_g \le \overline{p}_g^\text{gen},
\end{equation*}
with $\underline{p}_g^\text{gen}$ and $\overline{p}_g^\text{gen}$ both strictly positive. In our model of how flooding affects substations, this can easily lead to infeasibility. Suppose there exists a substation with one bus and that the bus has one generator and no load. If that bus is unaffected by flooding but all neighboring substations are affected, then the generator at the bus of the orphaned substation would have nowhere to send the strictly positive amount of power it is required to generate.

One way to eliminate this isolated generator condition is to employ an indicator variable $\gamma_g$ to track if generator $g$ is active or inactive and model
\begin{equation*}
    \underline{p}_g^\text{gen} \gamma_g \le \hat{p}_g \le \overline{p}_g^\text{gen} \gamma_g
\end{equation*}
such that the generator injects exactly zero power when it is inactive and must operate within its normal bounds when it is active. This is a reasonable model; however, it requires adding one discrete variable per generator to a second-stage power flow model that may otherwise have only continuous decision variables thus significantly increasing the complexity of the overarching two-stage model.

To avoid infeasibility in a computationally tractable way, we require each generator to produce active power within its lower and upper bounds, but we allow some or all of that power to be discarded at the point of generation. This effectively lowers the generator's lower active power bound to zero. The discarded power or overgeneration at a generator $g$ is a continuous variable denoted $\check{p}_g$, and we penalize the overgeneration in the objective function to incentivize solutions that do not involve overgeneration.

One interpretation of overgeneration is that power is arbitrarily discarded. This, of course, is a violation of energy conservation. Another interpretation of overgeneration is that generators are permitted to operate in a range that is likely to stress the hardware, potentially to the point of failure. We concede that overgeneration is not a real-world recourse option. However, and as we discuss further in Section~\ref{section:results}, we did not typically observe overgeneration in the solutions prescribed by our adapted models. Additionally, our evaluations of the prescribed mitigation solutions in an exact AC optimal power flow model (with overgeneration disallowed) demonstrate that the prescribed mitigation solutions are reasonable.

In practice, there are numerous controls for responding to severe contingencies. See, for example, the detailed ACPF-based contingency analysis model in \cite{Coffrin2019}. However, incorporating so many controls, especially discrete topology-affecting controls, in the recourse problem of a two-stage model may lead to computational intractability. For our time-sensitive application, we avoid this issue by rather incorporating a single indicator variable $\chi$. The purpose of $\chi$ is to trivially ensure feasibility when the fixed-topology recourse problem, fixed at least for given $\bs{x}$ and $\bs{\xi}$, is otherwise infeasible. Essentially, $\chi$ is a computationally tractable mechanism for ensuring relatively complete recourse by simply indicating that recourse options not captured by our model must be exercised to attain a feasible power flow solution. If not for relatively complete recourse, a solver would fail to prescribe a mitigation solution if the recourse problem for one or more scenarios were infeasible. After the adoption of a mitigation solution, of course, a more detailed model like the one in \citet{Coffrin2019} may be used to determine an effective response upon the realized flooding.

For the DC model, allowing load shed and overgeneration is sufficient for ensuring feasibility in well-parameterized instances. Regardless, we integrate $\chi$ in our adaptation of the DC model to demonstrate how to ensure feasibility if, for example, one elected to disallow overgeneration. For the more complex LPAC model, the existence of a non-trivial solution is not evident. We integrate $\chi$ in our adaptations so that the solution having $\chi = 1$, all $\alpha_n$ and $\beta_{nm}$ as dictated by constraints \eqref{eq:def_alpha} and \eqref{eq:def_beta}, all $\widehat{\cos}_{nm} = 1$, and all other variables equal to zero is feasible. We dub this the ``trivial solution'' as it features zero power generation, zero power flow, and zero load satisfied. There may be other feasible solutions having $\chi = 1$, but the trivial solution is always optimal in this case. From constraints \eqref{eq:infeasibility_load_shed}, $\chi = 1 \implies \delta_d = 0, \forall d \in D$. Ergo, the first term in the objective must be $\sum_{d \in D} p_d^\text{load}$ in any optimal solution involving $\chi = 1$. Additionally, having $\chi = 1$ reduces the effective lower bounds on active power generation to zero in constraints \eqref{eq:dc_generation_limits} from the DC model and constraints \eqref{eq:lpac_p_generation_limits} from the LPAC model. As such, there is no need for overgeneration, and the second term in the objective function must be zero in any optimal solution involving $\chi = 1$. Note that we incorporate $\chi$ in constraints \eqref{eq:lpac_p_ohms_law}-\eqref{eq:lpac_q_ohms_law} simply to ensure the feasibility of the trivial solution. In those constraints, having $\chi = 1$ nullifies the constant terms and the terms involving $\widehat{\cos}_{nm}$.

Recourse models \eqref{eq:L_dc} and \eqref{eq:L_lpac} are designed to admit a non-trivial PF solution with partial or full load satisfaction if such a solution exists and otherwise admit the trivial solution that features zero load satisfaction. This allows us to evaluate flooding scenarios of any severity without fear of the overarching two-stage model becoming infeasible. It also grants us the ability to effectively warmstart a solver with any feasible first-stage solution $\bs{x} \in \mathcal{X}$.

\section{Results} \label{section:results}

\newcommand{\pathresults}{results}

In this section, we describe the case studies we used and the experiments we executed to draw comparisons between the various adaptations of the DC and LPAC approximations that we defined in Section~\ref{section:modeling} as well as the \eqref{eq:SP} and \eqref{eq:RO} models. We then assess the optimal objective value and mitigation solution data we obtained.

\subsection{Case Studies} \label{subsection:casestudies}

We leverage two flooding case studies we developed in past work \citep{Austgen2023}. The first is based on Tropical Storm Imelda from 2019 and the second on Hurricane Harvey from 2017. These events both affected the Texas coastal region, especially Houston, and were noteworthy for their high volumes of precipitation \citep{Blake2018,Latto2020}. Both case studies incorporate an electrically equivalent reduced network based on the ACTIVS 2000-bus synthetic grid of Texas comprising 362 substations, 663 buses, and 1509 branches near the Texas coast. All parameters associated with flooding and the reduced grid are provided in the case studies. 
In these case studies, resilience level $r=1$ corresponds to a dam height of 0.534 meters and level $r=2$ to exactly 1 meter. The unattainable level of resilience is $\hat{r} = 3$ (\textit{i.e.}, flooding in excess of 1 meter is deemed inexorable). The marginal resource costs of mitigation are assumed to be increasing. Moreover, substations are partitioned into small, medium, and large categories based on the highest-voltage component. Marginal resource costs (\textit{i.e.}, each $c_{kr}$) for various substation sizes and levels of mitigation are summarized in Table~\ref{tab:mitigation_costs}.
\begin{table}
\TABLE
{Marginal Resource Requirements for Protecting Substations of Different Sizes \label{tab:mitigation_costs}}
{
    \begin{tabular}{|r|ccc|}
        \hline
        Highest-Voltage & \multicolumn{3}{c|}{Resilience Level} \\
        \cline{2-4}
        Component & $r=1$ & $r=2$ & $r=3$ \\
        \hline
        115 kV or 161 kV & 1 & 2 & 3 \\
                  230 kV & 2 & 4 & 6 \\
                  500 kV & 3 & 6 & 9 \\
        \hline
    \end{tabular}
}
{}
\end{table}

Additionally, the Imelda and Harvey case studies respectively comprise 4 and 25 equiprobable scenarios. The flood level indicators for these scenarios as well as the EV and MV scenarios are illustrated in Appendix C in the online supplement.
\subsection{Experiments and Budget Thresholds} \label{subsection:experiments}

We applied each two-stage model from Section~\ref{subsection:two_stage_models} to the Tropical Storm Imelda and Hurricane Harvey case studies while varying the resource budget $f$ to be integer values between zero and the model's budget threshold (\textit{i.e.}, the maximum useful budget) from Table~\ref{tab:study_extent}. For the Tropical Storm Imelda case study, we performed these experiments incorporating each of the adapted DC, LPAC-C, LPAC-F, and QPAC models defined in Section~\ref{subsection:power_flow_modeling}. The Hurricane Harvey case study comprises a greater number of scenarios that affect a larger number of substations as compared to the Imelda case study and thus presented a more difficult computational challenge. To compensate, we only incorporated the adapted DC and LPAC-C models to the Harvey-based instances.
\begin{table}
\TABLE
{Summary of Two-stage Models, Power Flow Models, Case Studies, and Budget Thresholds in our Experiments \label{tab:study_extent}}
{
    \begin{tabular}{|r|c|c|r|r|}
        \cline{2-5}
        \multicolumn{1}{r|}{} & \multicolumn{2}{c|}{Model} & \multicolumn{2}{c|}{Case Study} \\
        \cline{2-5}
        \multicolumn{1}{r|}{} & Name & Formulation / Description & Imelda & Harvey \\
        \hhline{|=|=|=|=|=|}
        \multirow{6}{*}{\rotatebox[origin=c]{90}{\parbox{2.9cm}{Two-Stage Model}}}
        & EEV & $\vphantom{\sum_{\substack{0\\0}}^{0}} \sum_{\omega \in \Omega} \prob(\omega) \mathcal{L}(\overline{\bs{x}}, \bs{\xi}^\omega)$ & 9 &  66 \\
        \cline{2-5}
        &  SP & $\vphantom{\sum_{\substack{0\\0}}^{0}} \min_{\bs{x} \in \mathcal{X}} \sum_{\omega \in \Omega} \prob(\omega) \mathcal{L}(\bs{x}, \bs{\xi}^\omega)$ & 20 & 193 \\
        \cline{2-5}
        & EWS & $\vphantom{\sum_{\substack{0\\0}}^{0}} \sum_{\omega \in \Omega} \prob(\omega) \min_{\bs{x} \in \mathcal{X}} \mathcal{L}(\bs{x}, \bs{\xi}^\omega)$ & 11 &  66 \\
        \cline{2-5}
        & MMV & $\vphantom{\min_{\substack{0\\0}}^{0}} \max_{\omega \in \Omega} \mathcal{L}(\widehat{\bs{x}}, \bs{\xi}^\omega)$ &  8 &  62 \\
        \cline{2-5}
        &  RO & $\vphantom{\min_{\substack{0\\0}}^{0}} \min_{\bs{x} \in \mathcal{X}} \max_{\omega \in \Omega} \mathcal{L}(\bs{x}, \bs{\xi}^\omega)$ & *9 & *62 \\
        \cline{2-5}
        & MWS & $\vphantom{\min_{\substack{0\\0}}^{0}} \max_{\omega \in \Omega} \min_{\bs{x} \in \mathcal{X}} \mathcal{L}(\bs{x}, \bs{\xi}^\omega)$ & *5 & *48 \\
        \hhline{|=|=|=|=|=|}
        \multirow{6}{*}{\rotatebox[origin=c]{90}{\parbox{2.90cm}{Power Flow Model}}}
        &     DC & \parbox{8.25cm}{\centering ~\\[-0.3em]
                                   sine is linear, cosine is 1, no reactive power \\
                                   branch conductance is negligible \\[0.2em]} & \checkmark & \checkmark \\
        \cline{2-5}
        & LPAC-C & \parbox{8.25cm}{\centering ~\\[-0.3em]
                                   sine is linear, cosine is 1, \\
                                   4-sided polyhedral relaxation of unit circle \\[0.2em]} & \checkmark & \checkmark \\
        \cline{2-5}
        & LPAC-F & \parbox{8.25cm}{\centering ~\\[-0.3em]
                                   sine is linear, 8-sided polyhedral relaxation of cosine, \\
                                   12-sided polyhedral relaxation of unit circle \\[0.2em]} & \checkmark & \\
        \cline{2-5}
        &   QPAC & \parbox{8.25cm}{\centering ~\\[-0.3em]
                                   sine is linear, quadratic relaxation of cosine, \\
                                   exact quadratic representation of unit circle \\[0.2em]} & \checkmark & \\
        \hline
    \end{tabular}
}
{\emph{Notes}: Values marked with an asterisk (*) were determined via optimization whereas those without were precomputed. In the row for \eqref{eq:EEV}, $\overline{\bs{x}}$ denotes the EV solution we identified. Likewise in the row for \eqref{eq:MMV}, $\widehat{\bs{x}}$ denotes the MV solution we identified.}
\end{table}
We implemented the extensive form (\textit{i.e.}, deterministic equivalent) of all the models from Table~\ref{tab:study_extent} in Python using the \texttt{gurobipy} package. To solve, we used the Gurobi solver \citep{Gurobi2022}, configured to terminate at a 0\% MIP gap. For the SP instances, we leveraged our greedy heuristic from \citet{Austgen2023} to warmstart the solver with mitigation solutions. The heuristic iteratively commits to the mitigation decision that restores the largest expected weighted combination of transmission capacity and load per \tigerdam~resource, and it terminates when no cost-feasible mitigation options remain. For the RO instances, we adapted the heuristic so that it instead restores the largest weighted combination of transmission capacity and load per resource in the approximate worst-case scenario. At each iteration, the approximate worst-case scenario is based on the difference between available generation capacity and serviceable load under the already-committed mitigation decisions. To further aid the solver, one could fix $x_{kr} = 0$ if $\xi_{kr}^\omega = 0$ for all $\omega \in \Omega$ (\textit{i.e.}, one could prohibit excessive substation flood mitigation). Even without such variable fixing, the Gurobi solver identifies these cases and eliminates the corresponding discrete variables in preprocessing.

We evaluated all the identified optimal mitigation solutions to the SP and RO instances in an ACPF model to gauge the efficacy of our approach. For this task, we adapted the polar power voltage (PSV) implementation of ACPF from the Python package EGRET \citep{Knueven2020} to our application. Specifically, we incorporated load shed variables at each bus, we modified the objective function so that it minimizes total active power load shed, and we imposed that active and reactive power are shed in the same proportion, which in our models is ensured by $\delta_d$. In the adapted ACPF model, overgeneration is strictly prohibited. To evaluate each mitigation solution in the ACPF model, we took the following steps:
\begin{enumerate}
    \item We obtained an initial power flow solution for each scenario by fixing the mitigation solution in the model that prescribed it and solving the independent recourse problems.
    \item On a scenario-by-scenario basis, we eliminated generators isolated from loads by the contingency to avoid the source of infeasibility described in Section~\ref{subsection:relatively_complete_recourse}.
    \item For each scenario, we identified an ACPF solution using IPOPT, a solver designed to find locally optimal solutions to continuous nonconvex optimization problems \citep{Wachter2006}.
\end{enumerate}

We performed all experiments on an SKX compute node from the Stampede2 cluster at the Texas Advanced Computing Center. The SKX compute node features dual Intel Xeon Platinum 8160 CPUs (24 cores, 2.10 GHz) and 192 GB of DDR4 RAM \citep{Stampede2}. We have uploaded the input data, code, and result data to a public GitHub repository hosted by the INFORMS Journal on Computing \citep{Austgen2024code}.

The resource budget thresholds defined in Table~\ref{tab:study_extent} are the number of resources that may be effectively used to prevent flooding in each of the two-stage models. Most of these values may be determined prior to solving the associated two-stage model simply by assessing the set of flooding scenarios and the resource costs of the possible mitigation actions. The two exceptions are the budgets for the RO model and MWS bound since the worst-case scenario is not known \textit{a priori} but is rather a function of mitigation decision making such that the threshold may only be determined through optimization. These values are preceded by an asterisk in Table~\ref{tab:study_extent}. While the thresholds for the RO model and MWS bound could differ based on the chosen PF model, we found them to be the same. To reduce verbosity, we use $\overline{f}_\text{X}$ throughout this section to denote the budget threshold for model X.

Recall that some flooding in our model is deemed inexorable via constraints \eqref{eq:con_inexorable}. The flooding in the MV scenario is necessarily worse than the flooding in the EV scenario, so the MV scenario has at least as much inexorable flooding. Because the threshold is based on resources used to effectively prevent flooding, it is possible for MMV bound's budget threshold to be less than the EEV bound's threshold. Additionally, because these thresholds are easily computed prior to optimization, we only had to compute the EV and MV solutions for integer budgets between 0 and the corresponding threshold. For budgets above the threshold, we avoided the arbitrary deployment of mitigation resources by simply adopting the solution corresponding to the budget threshold. The practical implication of our handling the budget threshold in this way is that the EEV and MMV bounds stagnate for budgets above the threshold. We also adopted this treatment of the budget threshold for the independent subproblems that must be solved to compute the EWS and MWS bounds.
\subsection{Objective Value and Bounds}

\begin{figure}
\FIGURE
{
    \ifdefined\epsswitch
        \includegraphics[width=6.5in]{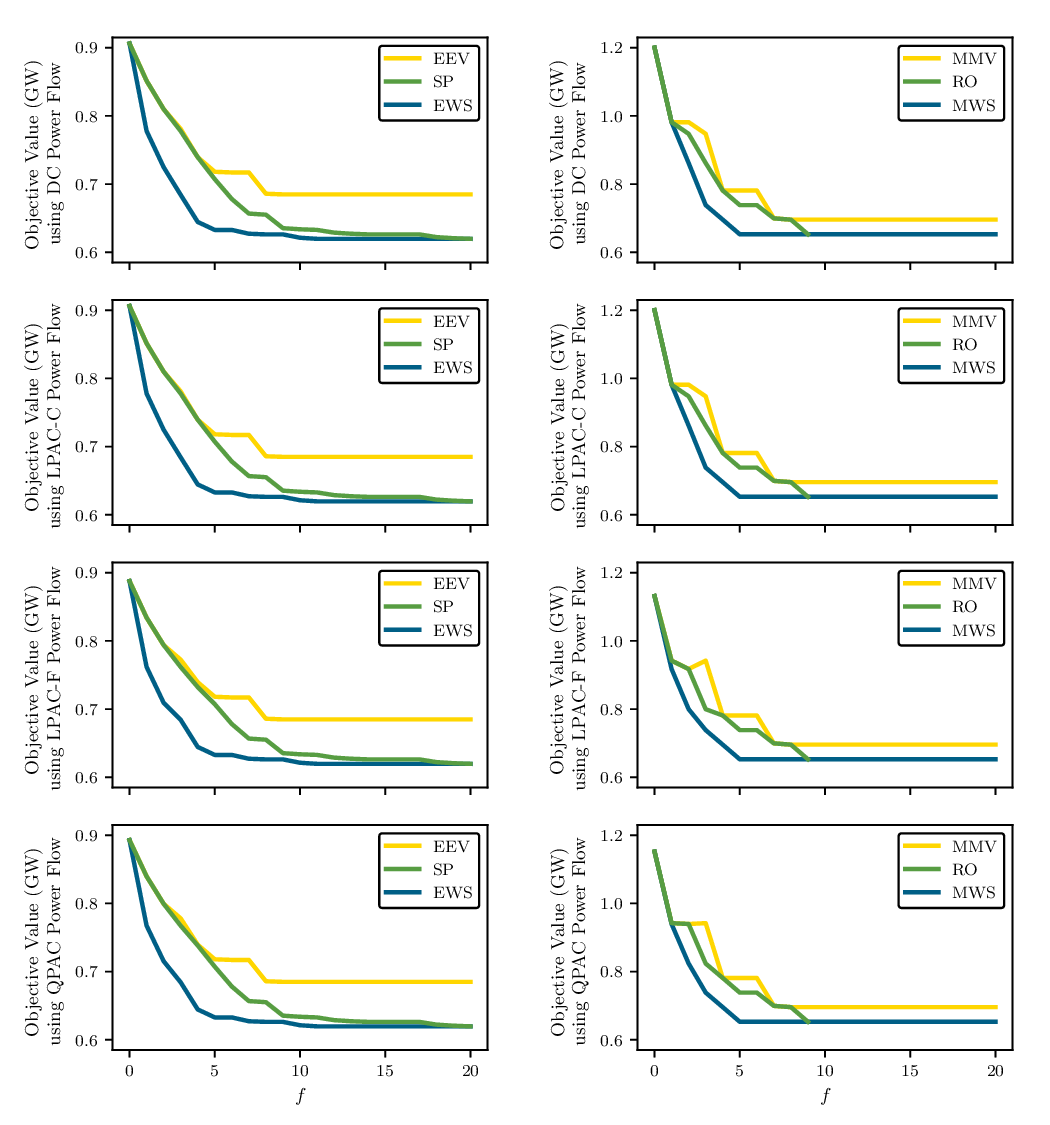}
    \else
        \includegraphics[width=6.5in]{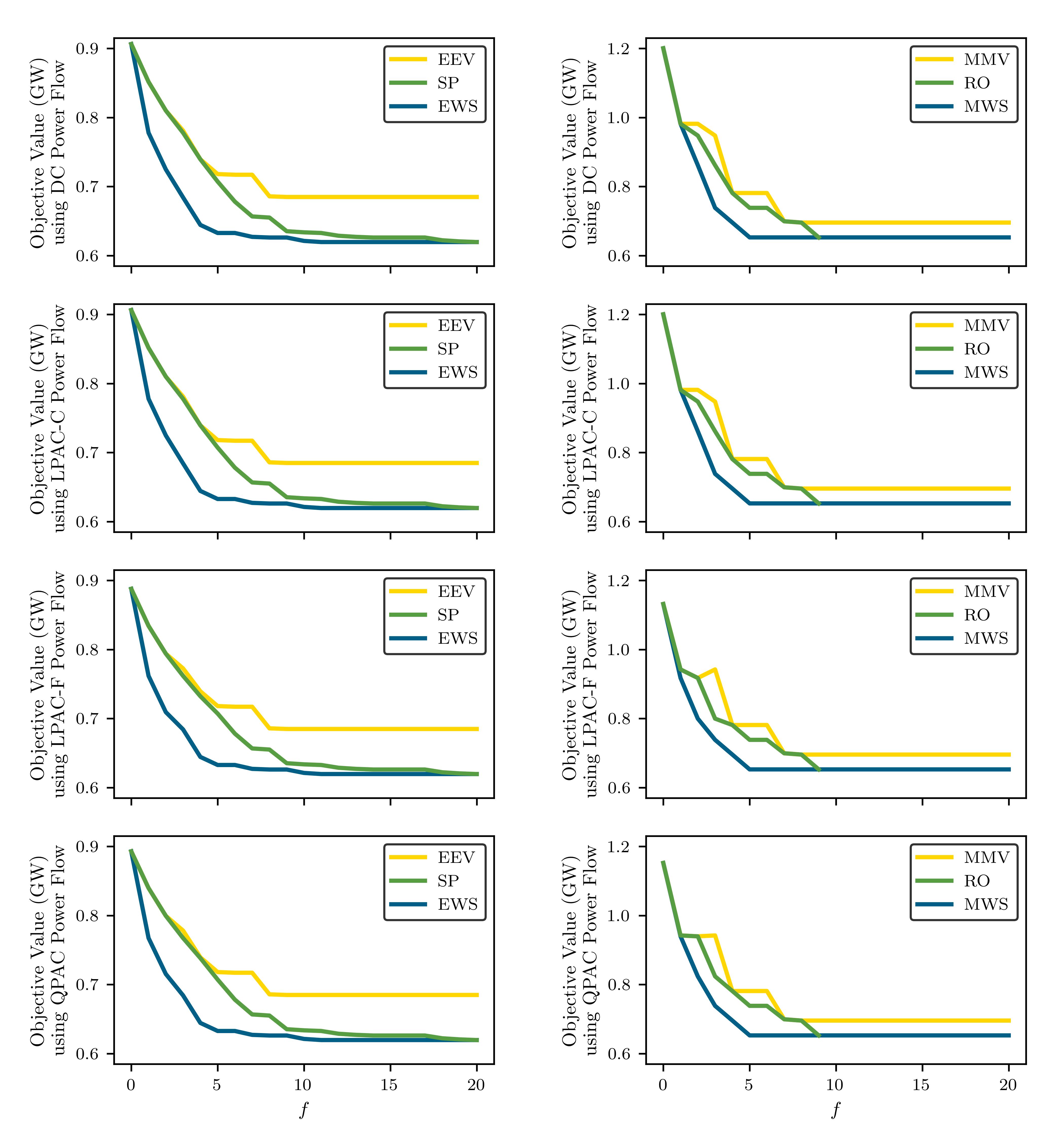}
    \fi
}
{Objective Values and Bounds as Functions of the Mitigation Budget in the Tropical Storm Imelda Case Study \label{fig:imelda_obj_bnds}}
{}  
\end{figure}

\begin{figure}
\FIGURE
{
    \ifdefined\epsswitch
        \includegraphics[width=6.5in]{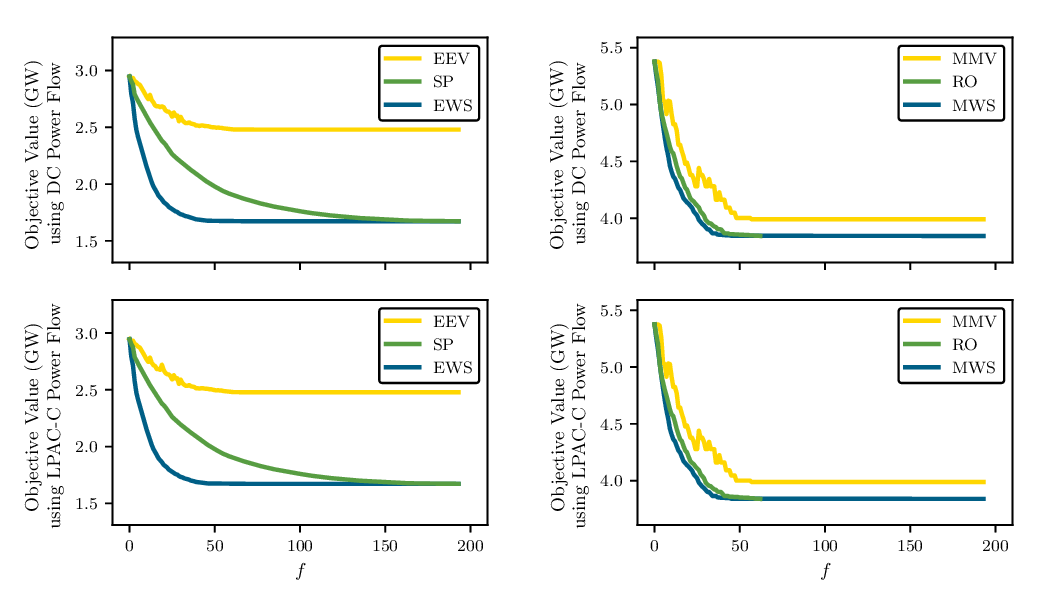}
    \else
        \includegraphics[width=6.5in]{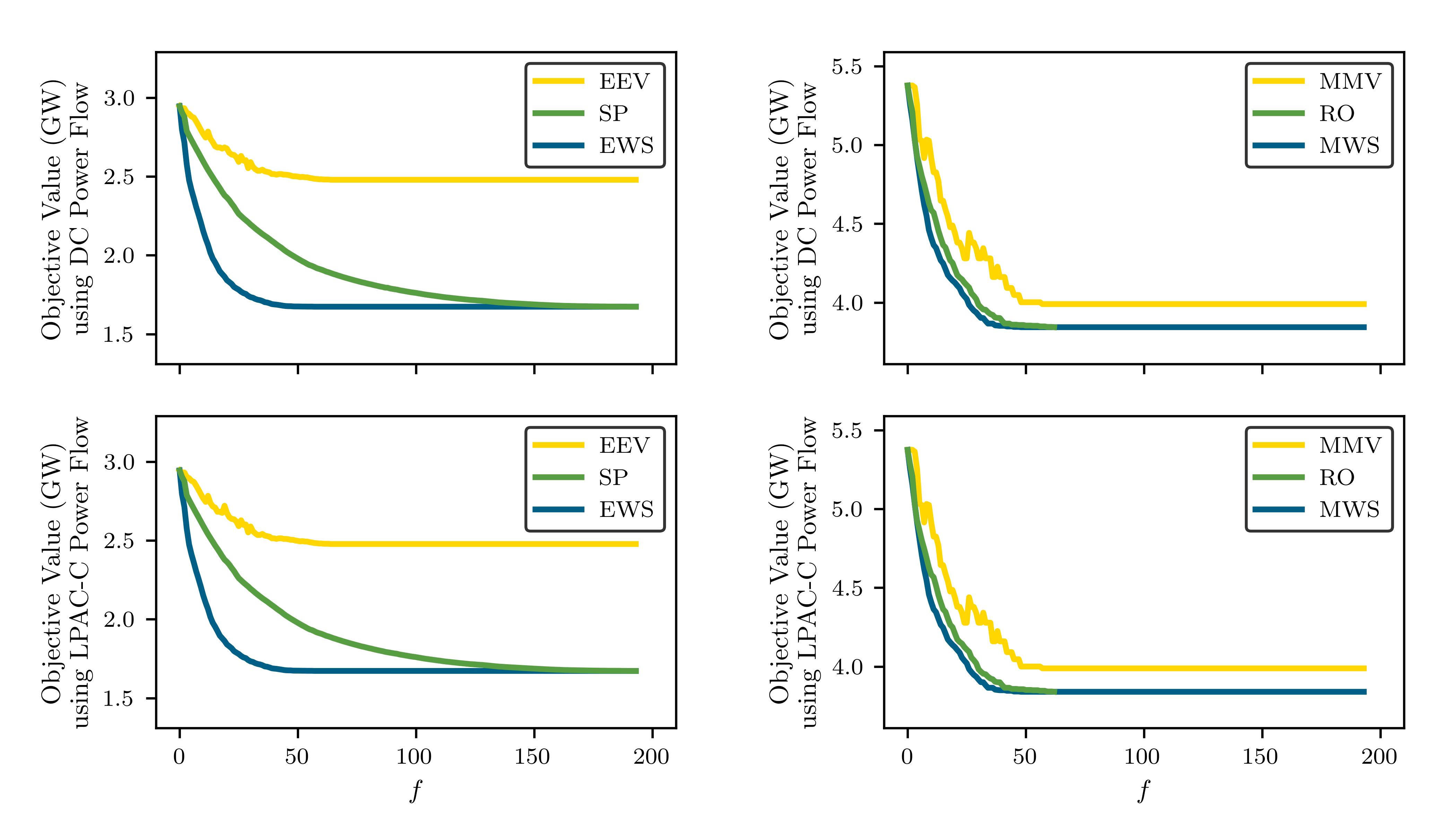}
    \fi
}
{Objective Values and Bounds as Functions of the Mitigation Budget in the Hurricane Harvey Case Study \label{fig:harvey_obj_bnds}}
{}
\end{figure}

We study the objective values of the two-stage models mainly to assess the trade-off between the resource budget and the achievable load shed. The objective values are also useful for comparing the effect of using different PF models, evaluating the tightness of the bounds on the SP and RO models, and computing quantities of theoretical interest like VSS and EVPI. Figures~\ref{fig:imelda_obj_bnds} and \ref{fig:harvey_obj_bnds} show the objective value data for the Imelda and Harvey case studies, respectively. Throughout the rest of this paper, we let $z_\text{X}^*$ denote the optimal objective value of model X and let $x_\text{X}^*$ similarly denote a mitigation solution capable of attaining that objective value. These quantities, of course, also depend on the other experimental controls like the chosen power flow model and the budget $f$. In our discussion of these quantities, we make their dependence on other controls clear via prose and context.

As expected, $z_\text{SP}^*$ and $z_\text{RO}^*$ monotonically decrease in the budget $f$. This is a theoretical guarantee since increasing $f$ enlarges the solution space thus allowing better mitigation strategies. For the same reason,  $z_\text{EWS}^*$ and $z_\text{MWS}^*$ are also monotonic. Because $z_\text{EEV}^*$ and $z_\text{MMV}^*$ are based on evaluations of fixed first-stage solutions, these bounds are not monotonically decreasing in $f$, though they are generally decreasing. Both $z_\text{SP}^*$ and $z_\text{RO}^*$ initially make steady improvements as the budget increases from zero, though the marginal improvements granted by additional resources are generally decreasing. In both case studies, $z_\text{EEV}^*$ and $z_\text{EWS}^*$ are poor bounds on $z_\text{SP}^*$ compared to how well $z_\text{MMV}^*$ and $z_\text{MWS}^*$ bound $z_\text{RO}^*$.

For Imelda, Figure~\ref{fig:imelda_obj_bnds} shows that $\overline{f}_\text{SP} = 20$ \tigerdams~stacked up to one meter tall are able to reduce expected loss (\textit{i.e.}, $z_\text{SP}^*$) by 31.6\% compared to the ``do-nothing'' solution. The worst-case loss (\textit{i.e.}, $z_\text{RO}^*$) may be reduced by 42.3\% when $\overline{f}_\text{RO} = 9$ such \tigerdams~are deployed. For Harvey, Figure~\ref{fig:harvey_obj_bnds} shows improvements of 43.2\% and 28.5\%, respectively, when $\overline{f}_\text{SP} = 193$ or $\overline{f}_\text{RO} = 62$ \tigerdams~are deployed. In these figures, note that the objective value decreases to a strictly positive value rather than to zero as a consequence of the case studies featuring inexorable flooding.

The bounds in the Imelda case study are tighter than those in the Harvey case study. This is reasonably explained by each $\bs{\xi}^\omega$ for $\omega \in \Omega$ being similar to the others and also to $\overline{\bs{\xi}}$ and $\hat{\bs{\xi}}$. This leads to the first-stage decisions in the EV solution and optimal solutions to the SP model and EWS bound subproblems all being similar to one another, and the same is true for the MV solution and optimal solutions to the RO model and MWS bound subproblems.

In both case studies, the VSS (\textit{i.e.}, the distance between the yellow $z_\text{EEV}^*$ and green $z_\text{SP}^*$ series) is generally increasing in the budget and largest when the budget is near $\overline{f}_\text{SP}$. That is, the SP model is more valuable for planning with larger budgets than for smaller budgets. The size of the gap also indicates that planning for the expected flooding scenario is a poor strategy, especially for larger budgets. In contrast, the EVPI (\textit{i.e.}, the distance between the green $z_\text{SP}^*$ and blue $z_\text{EWS}^*$ series) rapidly increases from zero then steadily decreases to zero as the budget increases. Specifically, EVPI is largest for budgets roughly halfway between zero and $\overline{f}_\text{EWS}$. Intuitively, budgets near the halfway mark are small enough that the marginal return on investment is still substantial and large enough that the variation in flooding scenarios may be more effectively addressed by nuanced decision making. For the RO model, $z_\text{MMV}^* - z_\text{RO}^*$ is largest for budgets roughly halfway between zero, and $\overline{f}_\text{MMV}$, and $z_\text{RO}^* - z_\text{MWS}^*$ is largest roughly halfway between zero and $\overline{f}_\text{MWS}$. Though the reason for the former phenomenon is difficult to determine, the reason for the latter is conceivably the same as for the EVPI.

\begin{figure}
\FIGURE
{
    \ifdefined\epsswitch
        \includegraphics[width=6.5in]{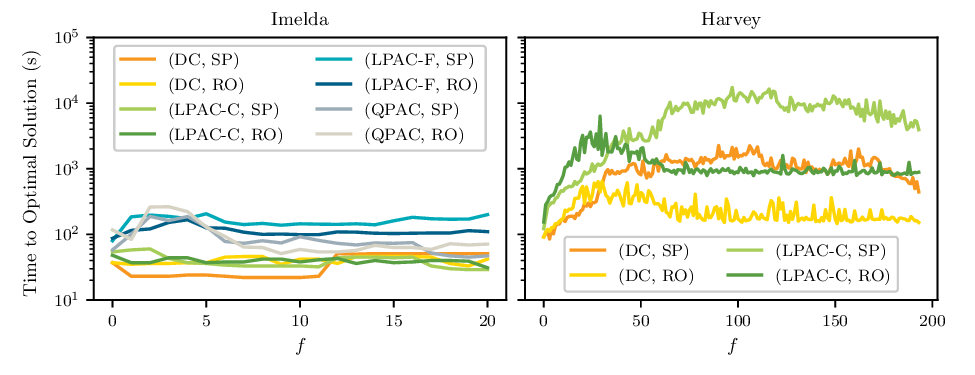}
    \else
        \includegraphics[width=6.5in]{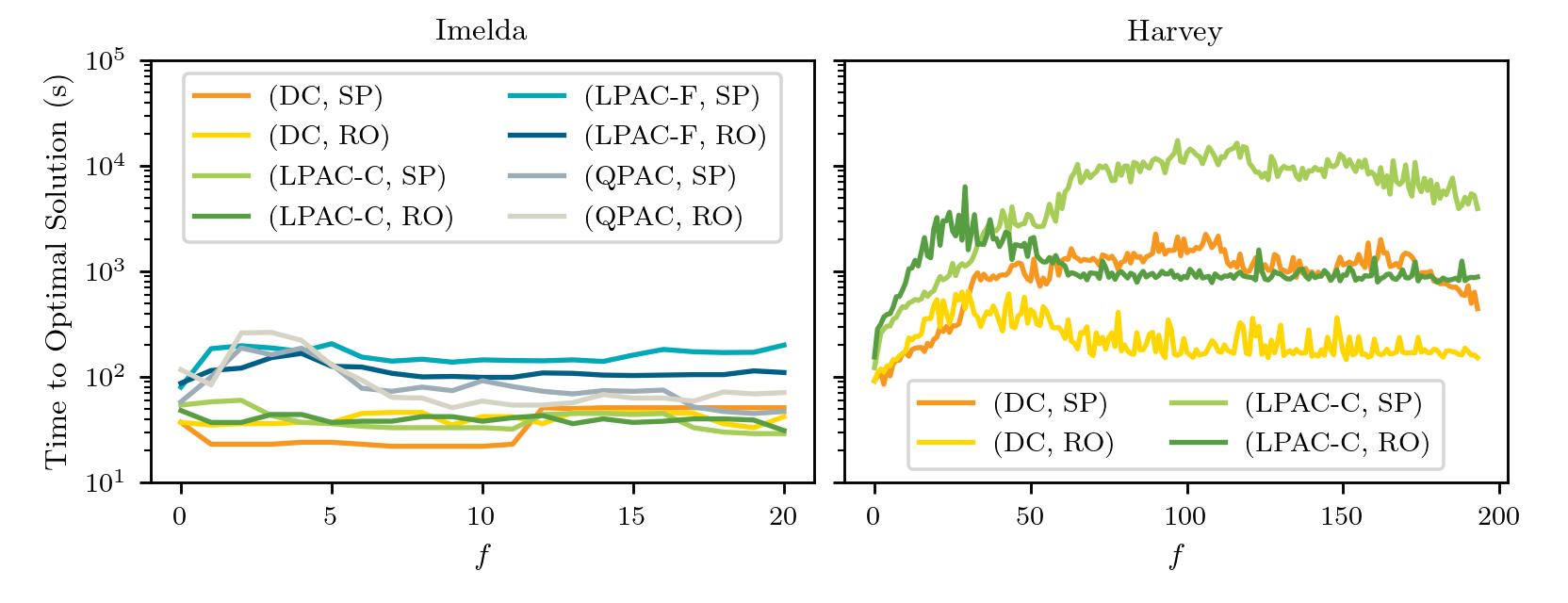}
    \fi
}
{Times Required to Solve the Studied Instances to Optimality \label{fig:solution_times}}
{}
\end{figure}

For the Imelda case study, the time required by the solver to identify a solution and prove its optimality, henceforth referred to as the ``time to optimal solution,'' was not affected much by the resource budget $f$ as indicated in Figure~\ref{fig:solution_times}, and all instances solve in a short time. The time to optimal solution for Harvey-based instances, however, was more dependent on $f$. For example, the times for the SP model with LPAC-C increased from roughly 2 minutes for $f = 0$ to nearly 5 hours for $f = 97$. Independent of the incorporated PF model, the instances of the SP and RO model that are fastest to solve are those with a small budget, which makes sense given that the size of the mitigation solution space is strictly increasing in the budget via constraint \eqref{eq:con_resource_hi}. Interestingly, however, the instances of the SP model that are slowest to solve are not those with budgets near $\overline{f}_\text{SP}$, but rather those with budgets in the vicinity of $\overline{f}_\text{SP} / 2$. Similarly for the RO model, the instances that are slowest to solve are those with budgets in the vicinity of $\overline{f}_\text{RO} / 2$.

Differences induced by the various power models are subtle as shown in Figures~\ref{fig:imelda_obj_bnds} and \ref{fig:harvey_obj_bnds}. We attribute this to the fact that the PF models all share the same core constraints (\textit{i.e.}, KCL and an approximation of Ohm's Law) and are only subtly different themselves. No matter, while properties like the VSS and EVPI and the effect of the budget on the objective value are of direct interest, objective value differences induced by the chosen PF model only indicate that we should assess the mitigation solutions more qualitatively.
\subsection{Mitigation Decisions}

The mitigation decisions are the most important decisions in our model and the other decisions serve mainly to capture the operational recourse that influences the mitigation decisions. The operational decisions are for a proxy of ACPF, and the realized flooding is unlikely to exactly match any of the sampled scenarios anyway. Nevertheless, if the choice of surrogate PF model used in the recourse stage affects the optimal mitigation, then it raises the question if one or more of those models are not suitable for this application.

For the vast majority of the instances we studied, the optimal mitigation solutions were uniquely optimal. We determined this computationally by adding the no-good cut ${(\bs{1} - \bs{x}^*)^\top \bs{x} \ge 1}$ to the instance, solving that restricted instance to optimality, and assessing the objective value. In this cut, $\bs{x}^*$ is the identified optimal solution, $\bs{x}$ is the vector of decision variables, and $\boldsymbol{1}$ is the all-ones vector. For the vector of resource costs $\bs{c}$, if $\bs{c}^\top \bs{x}^* = f$, then the only solution cut from the feasible space is $\bs{x}^*$. Otherwise, this constraint cuts all $\bs{x} \ge \bs{x}^*$ (\textit{i.e.}, all the solutions that implement $\bs{x}^*$ and possibly more). Despite mitigating more flooding, such solutions are no better performing than $\bs{x}^*$, which is optimal.

We attribute the frequent occurrence of unique optimal mitigation solutions to the power grid instance having little symmetry. The few instances with multiple mitigation optima typically have one or more resources that cannot be used effectively and are instead deployed arbitrarily. As an example, consider a deterministic instance in which only one bus, a load bus, is affected. If 2 resources are required to prevent flooding at the associated substation, but only 1 resource is available, then that resource may be deployed arbitrarily with no effect. For the few instances having multiple mitigation optima, we did not bother identifying all optima as the task is computationally burdensome. Nevertheless, because most instances have a unique mitigation optimum, the vast majority of our mitigation solution comparisons are not subject to the doubt that arises from not having identified all of the multiple optima.

\subsubsection{Impact of Power Flow Model.}
To determine the impact of the choice of power flow model on the optimal mitigation decisions, we performed several pairwise comparisons of mitigation solutions prescribed under our various PF models. For the Imelda case study, we compared every pair of PF models for budgets from 1 to 20 in the SP model and from 1 to 9 in the RO model. For Harvey, we compared DC and LPAC-C for budgets from 1 to 193 in the SP model and from 1 to 62 in the RO model. This amounted to 429 comparisons. To assess the similarity of two mitigation solutions, we employed an absolute measure from our past work \citep{Austgen2021} and an analogous relative measure:
\begin{gather*}
    \abssim(\bs{x}_\text{A}, \bs{x}_\text{B}; \bs{c})
    = \bs{x}_\text{A}^\top \diag(\bs{c}) \bs{x}_\text{B}
    \label{eq:abssim}, \\
    \relsim(\bs{x}_\text{A}, \bs{x}_\text{B}; \bs{c})
    = \frac{\bs{x}_\text{A}^\top \diag(\bs{c}) \bs{x}_\text{B}}{\max\{\bs{c}^\top \bs{x}_\text{A}, \bs{c}^\top \bs{x}_\text{B}\}}.
    \label{eq:relsim}
\end{gather*}
In these definitions, $\diag(\cdot)$ denotes the matrix having the argument as its diagonal. The absolute measure, a weighted inner product, captures the number of resources deployed in the same manner. The relative measure normalizes the quantity to a proportion of the number of resources needed by the more resource-intensive solution. Among the 429 comparisons, only 29 had different solutions, which are detailed in Appendix D in the online supplement.

For solutions from instances with small resource budgets, the relative similarity tends to be smaller. This is expected, however, since a small absolute difference in resource allocation impacts the normalized similarity measure more when the solutions being compared use fewer resources. In contrast, the similarity of solutions from instances with large resource budgets tend to be much higher.

In cases where the mitigation solutions were similar but not the same, we hypothesized that similar mitigation solutions would perform similarly if fixed in a model regardless of which PF model were to be incorporated in the second stage. We tested this hypothesis for all 29 cases and found it to be true. In fact, we found differing solutions performed comparably even when they were not remarkably similar. For all pairs of differing mitigation solutions, we evaluated each solution in the model that prescribed the other and observed the optimality gaps to be quite small. The optimality gaps never exceeded 3.6\% for Imelda instances and never exceeded 0.11\% for Harvey instances.

This analysis indicates that the added fidelity offered by an LPAC model does not significantly impact the decision making compared with a DC model. Though this is a key finding of this paper, it is important to qualify given that in \citet{Coffrin2014} the LPAC model was found to admit more accurate PF solutions than the DC model.
First, it is important to note that in this application load may only be served at operational substations, and that mitigation decision making is largely driven by the need to spare load-serving substations. The importance of sparing a specific substation is based in part on the set of all substations spared and the ability of the system to facilitate the flow of power to those locations. The likelihood and severity of flooding at each substation are nevertheless crucial factors independent of the incorporated PF model, so the application itself is such that the optimal mitigation solutions induced by different PF models may be similar or the same.
Second, in Section~\ref{subsection:pf_model_comparisons} we mention a number of power system applications in which the impact of the incorporated PF model has been assessed. For most of the applications, comparisons were based on the similarity of decisions pertaining or tightly linked to the PF model (\textit{e.g.}, operating, grid posturing, and restoration scheduling decisions), and that the choice of PF model affects such decisions is expected. Only \citet{Garifi2022} likewise compared PF models on the basis of the mitigation investments they induce, and they concluded that the investments may differ based on the incorporated PF model. We offer three potential reasons for our differing conclusion. First, our recourse problem features a single time period and does not consider the restoration process whereas the recourse problem in \cite{Garifi2022} has multiple time periods linked by generator ramp constraints and restoration scheduling decisions that might be more limiting under one PF model than another. Second, mitigation decisions in our model are made under uncertainty of the threat whereas in \cite{Garifi2022} only a single deterministic contingency is considered. The risk measures we use to aggregate loss over the scenario-based representation of the uncertainty, especially the expected value risk measure, might help smooth any differences that arise in individual scenarios. Third, and potentially most importantly, the mitigation solutions we identify are truly optimal with respect to the chosen surrogate PF model. In \cite{Garifi2022}, the authors reported having configured the Gurobi solver to terminate within a 1\% MIP gap. In our experiments, we observed the existence of numerous, not necessarily all similar, near-optimal mitigation solutions. It was for this reason that we elected to solve all instances of our models to a 0\% MIP gap and moreover report the evaluated performance of dissimilar mitigation solutions under alternative PF models (\textit{i.e.}, our results in Appendix B in the online supplement). Given that the illustrated differences in the optimal mitigation solutions in \cite{Garifi2022} appear minor, it seems reasonable based on our observations that the DC and SOCP models compared in that paper might have induced the same mitigation solution had the instances been solved to a 0\% MIP gap.

\subsubsection{Impact of Uncertainty Perspective.}

\begin{figure}
\FIGURE
{
    \includegraphics[width=6.5in]{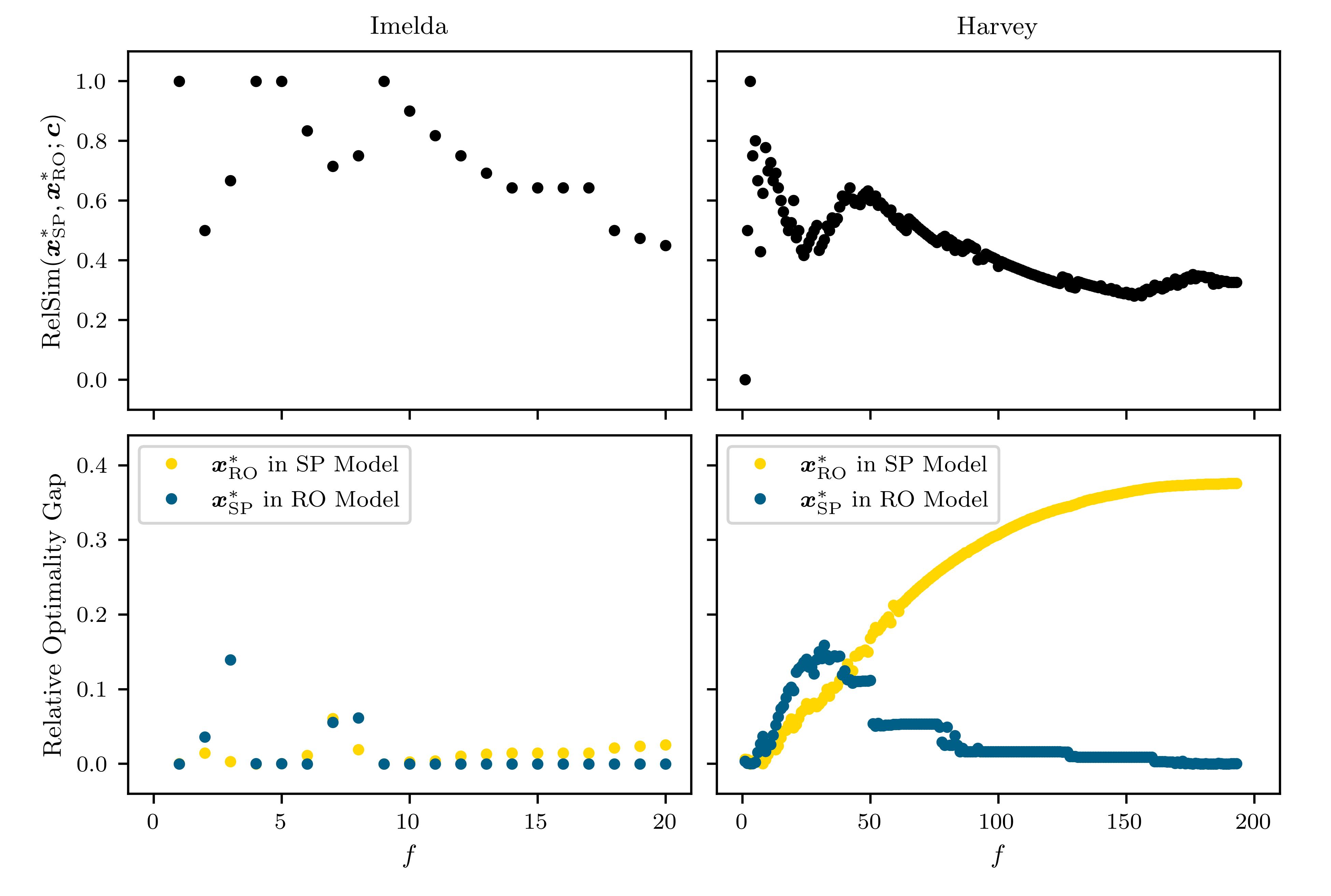}
}
{Similarity of the SP and RO Solutions and the Performance of Each in the Alternative Model
\label{fig:sp_ro_similarity}}
{}
\end{figure}

In Section~\ref{subsection:two_stage_models}, we discussed the uncertainty perspectives associated with the SP and RO models. We investigated how differently the two models, when subject to the same resource budget, prescribe mitigation decisions. Specifically, we quantified the similarity of the optimal mitigation solutions prescribed by the two models and moreover evaluated the performance of each prescribed solution under the alternative uncertainty perspective. For the models incorporating our adaptation of the DC approximation, these similarities are illustrated in Figure~\ref{fig:sp_ro_similarity}.

Because $\overline{f}_\text{RO} < \overline{f}_\text{SP}$ in both case studies, comparisons for $f > \overline{f}_\text{RO}$ involve the mitigation from the SP model with budget $f$ and the mitigation from the RO model with fixed budget $\overline{f}_\text{RO}$. Consequently, the relative optimality gap of $\bs{x}_\text{SP}^*$ in the RO model generally improves, and the relative optimality gap for $\bs{x}_\text{RO}^*$ in the SP model monotonically worsens as $f$ increases for $f > \overline{f}_\text{RO}$. In Figure~\ref{fig:sp_ro_similarity}, we see in the top two plots that optimal mitigation solutions of the two models often differ considerably with respect to how resources are allocated. In the bottom left plot, we observe in the Imelda case study that $\bs{x}_\text{RO}^*$ performs well in the SP model and likewise that $\bs{x}_\text{SP}^*$ performs well in the RO model. This is likely a consequence of the Imelda flooding scenarios being so similar. As we see in the bottom right plot, this phenomenon is not observed in the solutions from the Harvey case study, which has greater variance in the flooding scenarios. In the bottom plots of this figure, the relative optimality gap refers to the distance between the optimal objective value and the best attainable objective value when $\bs{x}$ is fixed to the optimal solution from the alternative model.

\begin{figure}
    \begin{subfigure}{0.49\textwidth}
        \centering
        \includegraphics[width=3in]{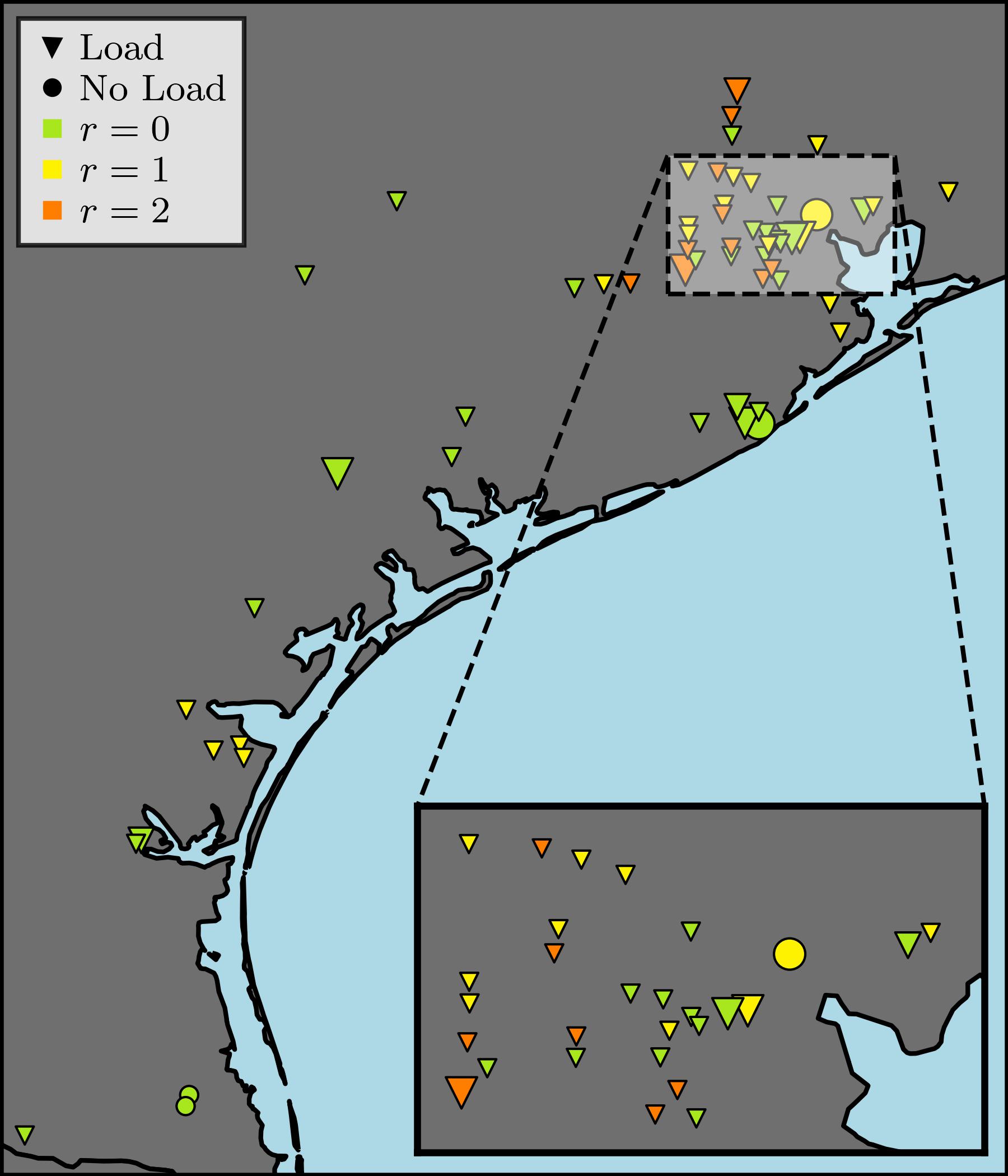}
        \caption{SP Solution}
        \label{fig:sp_map}
    \end{subfigure}
    \begin{subfigure}{0.49\textwidth}
        \centering
        \includegraphics[width=3in]{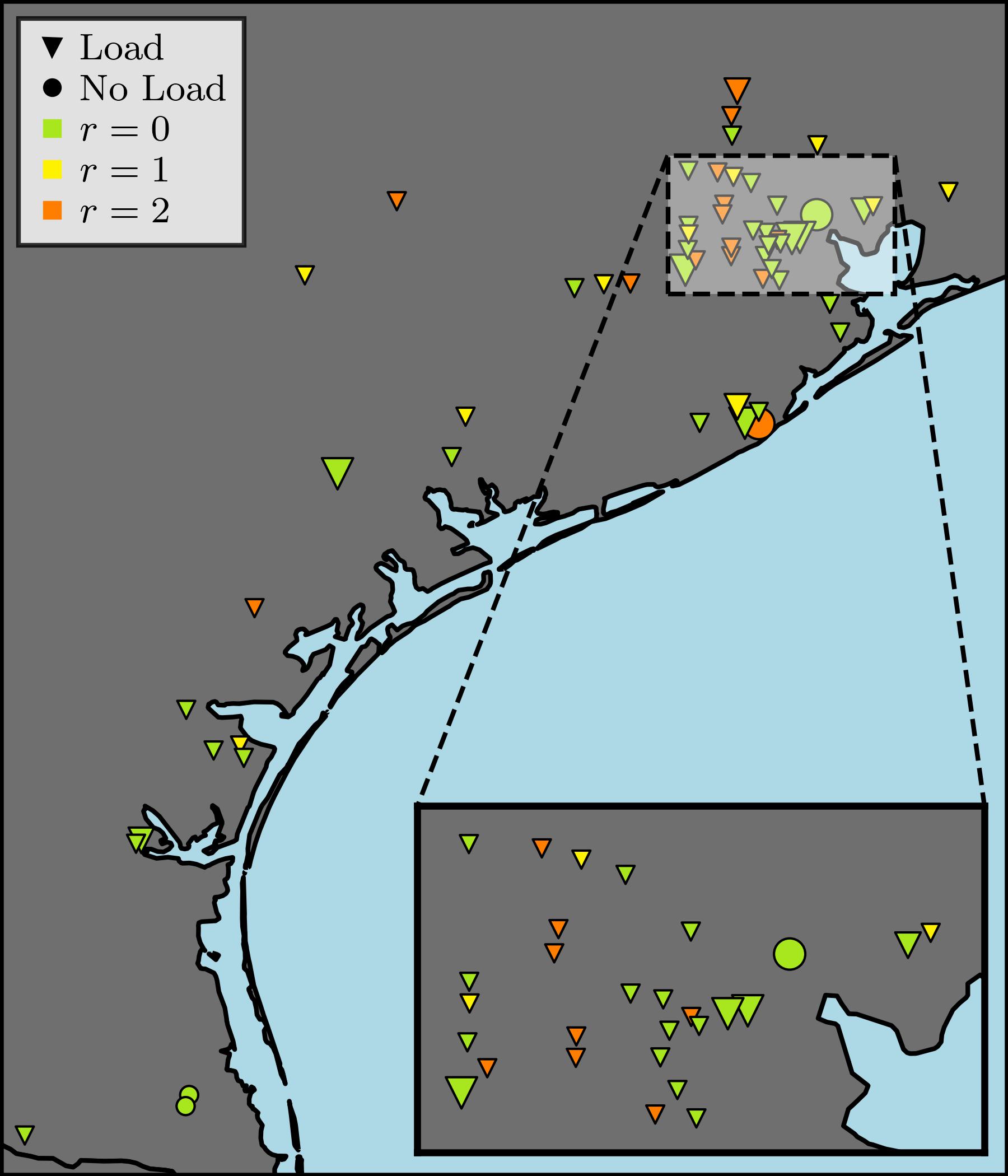}
        \caption{RO Solution}
        \label{fig:ro_map}
    \end{subfigure}
    \caption{Comparison of the SP and RO Mitigation Solutions for the Hurricane Harvey Instance with $f=62$ and $\hat{r} = 3$}
    \label{fig:sp_vs_ro_map}
\end{figure}

Qualitatively, $\bs{x}_\text{SP}^*$ is characterized by a relatively widespread allocation of the mitigation resources across the substations whereas $\bs{x}_\text{RO}^*$ tends to exhibit more concentrated allocation. In Figure~\ref{fig:sp_vs_ro_map}, this is illustrated for the case of $f = \overline{f}_\text{RO} = 62$ in the Harvey case study. On the left, we see $\bs{x}_\text{SP}^*$ sites level $r=1$ mitigation at 19 substations (yellow) and level $r=2$ at 10 (orange); on the right, however, $\bs{x}_\text{RO}^*$ sites level $r=1$ mitigation at only 10 substations and level $r=2$ at 14. Intuitively, the former solution is fit to minimize the expected loss over all flooding scenarios and thus must cover more assets. In contrast, the latter solution is fit to address only the scenarios that induce the worst outcomes.

\subsubsection{Evaluation in AC Power Flow Model}
\begin{figure}
\FIGURE
{
    \includegraphics[width=6.5in]{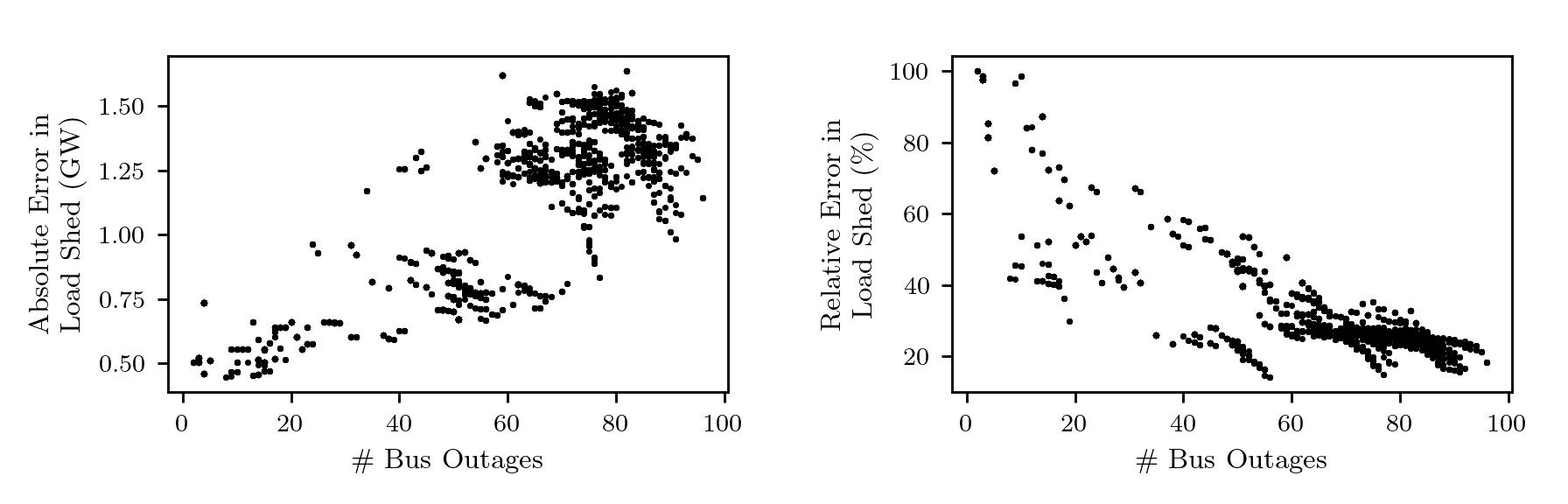}
}
{Absolute and Relative Error in Load Shed as a Function of Contingency Severity as Determined by AC Model Evaluations of the Optimal DC Model Solutions \label{fig:acpf-error}}
{}
\end{figure}

As mentioned in Section~\ref{subsection:experiments}, to validate our models we evaluated the performance of the mitigation solutions they prescribed under an AC power flow model. In these evaluations, we fixed the mitigation solutions and solved the independent recourse problems to local optimality. Taking the solutions of the AC power flow model as the ground truth, we observed that our models underestimate the total load shed by a substantial margin, similar to what was observed in \citet{Garifi2022}. As shown in Figure~\ref{fig:acpf-error}, we illustrate the absolute and relative error of the DC model.
In each of these two plots, each point represents a distinct contingency that arose under at least one of the optimal mitigation solutions our models prescribed.
As shown in the left plot of Figure~\ref{fig:acpf-error}, the absolute error in load shed ranged from 0.4 GW for minor contingencies to upwards of 1.6 GW for the most severe contingencies. Here, we measure the severity of the contingency by the number of bus outages. The absolute error was generally worse for more severe contingencies. Interestingly, as shown in the right plot of Figure~\ref{fig:acpf-error}, the relative error tended to be smaller for more severe contingency states. That is, though the error in total load shed increases absolutely in the severity of the contingency, it increases slowly relative to the increase in total load shed itself. Though we illustrate this here only for the DC model, we observed the same trends for LPAC-C, LPAC-F, and QPAC. In the AC solutions, we observed that most of the additional load shed is caused by bus voltage magnitude limits and branch thermal limits being reached in the vicinity of the buses where the additional load is shed. We conjecture that if these constraints are binding for certain buses and branches under normal circumstances, then substation flooding affecting those components may lead to the optimal AC power flow-informed mitigation solution being different from the solutions that are optimal under our surrogate models of power flow.

\begin{figure}
\FIGURE
{
    \includegraphics[width=6.5in]{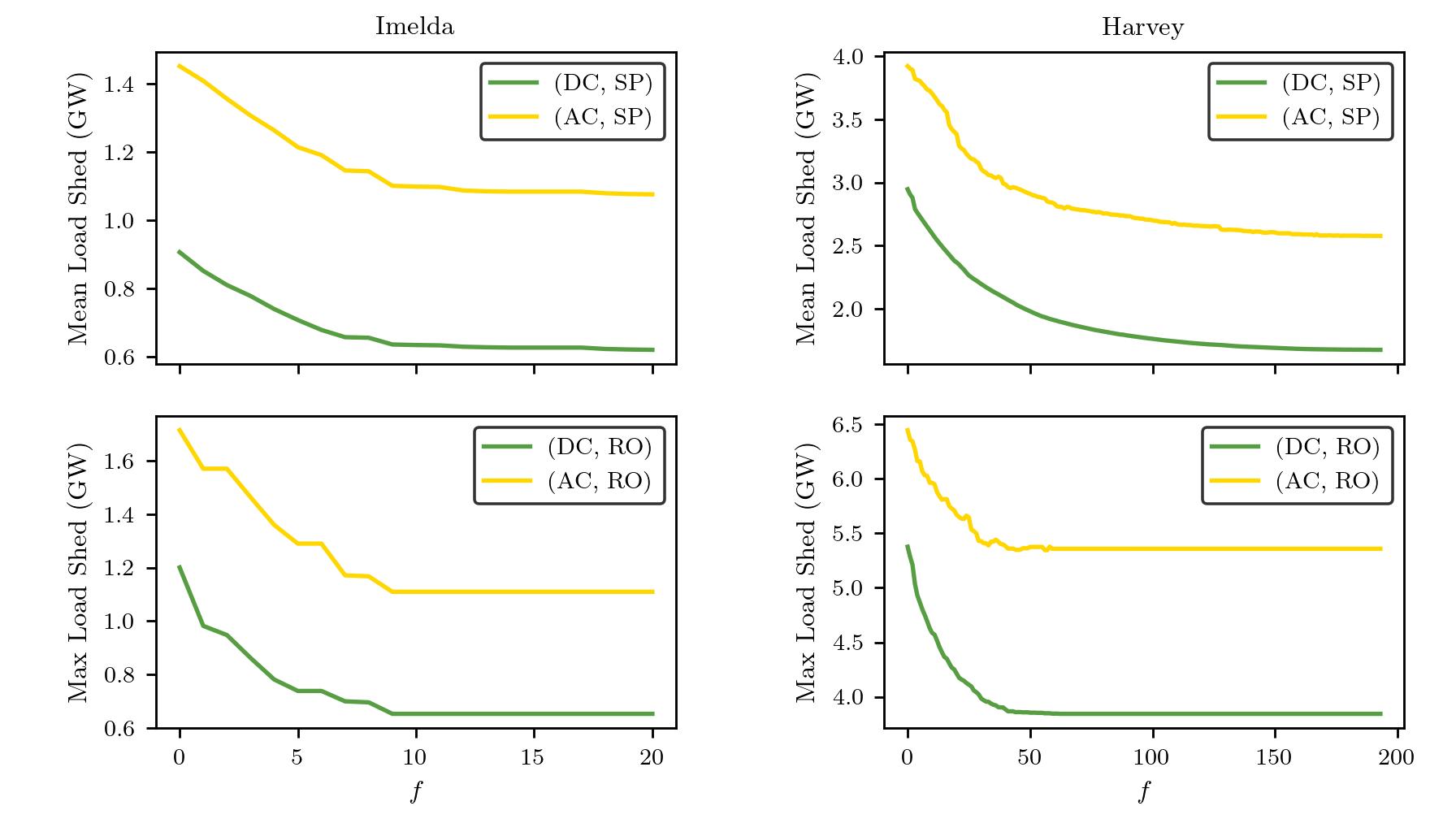}
}
{Load Shed Resulting From the Optimal DC Model Solutions Evaluated Using DC and AC Models \label{fig:acpf-perf}}
{}
\end{figure}

Despite our models incorporating surrogate PF models and overgeneration which lead to sizable approximation errors, they are evidently capable of prescribing effective mitigation solutions. For both case studies and uncertainty perspectives, we illustrate this for the DC model in Figure~\ref{fig:acpf-perf}.
The load shed in the AC model, as seen in this figure, was generally decreasing in the budget. The subtle nonmonotonicity in the load shed under the AC model could be a result of the mitigation being suboptimal or us having only obtained a local optimum to the AC-based recourse problem. Regardless, taking the ``no-mitigation'' load shed as the baseline, the mitigation solutions prescribed by our model improved the load shed as the mitigation budget increased and by roughly the same magnitude indicated by the AC model. Again, though we only illustrate this for the solutions prescribed by the DC model, we observed the same for LPAC-C, LPAC-F, and QPAC.

\section{Conclusions} \label{section:conclusions}

In this paper, we presented two-stage optimization models for informing substation flood mitigation decisions prior to an imminent hurricane. In a computational study, we applied our models to case studies featuring the ACTIVS 2000-bus synthetic grid of Texas and flooding scenarios derived from historical Tropical Storm Imelda and Hurricane Harvey data. We observed that compared to using a DC model, using an LPAC model has little to no bearing on the optimal mitigation decisions and objective values but has a remarkable impact on the time required to solve the model to optimality. Thus, we have reason to believe the DC model should be preferred to any of the LPAC models for this application despite its lower fidelity. We found the adopted uncertainty perspective (\textit{i.e.}, the decision maker's belief that either ``nature is fair'' or ``nature is adversarial'') significantly impacts the optimal mitigation solution. Also, we generally observed decreasing returns in the objective value as the mitigation budget increases. We demonstrated that our models, despite their inaccuracies, are capable of prescribing effective mitigation by evaluating the mitigation solutions using an AC power flow model. These results highlight the potential for \tigerdam~barriers to prevent power grid damages and for our models to prescribe their effective deployment.

Regarding future modeling efforts, we would like to extend the recourse problem to include multiple time periods, as it would allow us to more accurately capture temporal variability like the daily periodicity of demand and renewable generation and temporal constraints like generator ramp rate limits. This would make our models more difficult to solve, but more applicable to the challenges of modern power grids.
Also, we sought to improve resilience in this paper by minimizing system-wide load shed. Other considerations, such as the joint minimization of mitigation and load shed costs, equity-related metrics, chance constraints, etc., may be considered in the future.
Regarding future analysis efforts, understanding how the prescribed mitigation is affected by the scenario sampling scheme (\textit{e.g.}, by including or omitting stratification, adjusting the number of sampled scenarios, etc.) would help inform how the models may be most effectively used by practitioners. We are also interested in analyzing alternative power grid instances to understand if and how comparisons of the factors we assessed in this work are affected by specific network properties. Finally, we aim to incorporate convex relaxations of the AC power flow model, \textit{e.g.}, the semidefinite programming relaxation \citep{Molzahn2019}, to bound the performance of our models' prescribed solutions in the context of the AC power flow model.

\begin{appendices}

\section{Tightening the Polygonal Relaxation of Cosine}
\renewcommand{\pathmodeling}{modeling}
Let
\begin{equation*}
    \mathcal{B}_1(\theta; \hat{\theta}) = (\hat{\theta} - \theta) \sin(\hat{\theta}) + \cos(\hat{\theta})
\end{equation*}
be the line tangent to $\cos(\theta)$ at $\hat{\theta}$.
In \citet{Coffrin2014}, the LPAC model incorporates
\begin{equation*}
    \widehat{\text{cos}}_{nm} \le \mathcal{B}_1(\theta_n - \theta_m; \hat{\theta}),
    \quad \forall \hat{\theta} \in \Theta_\text{cos} = \left\{\left(\tfrac{2t - T - 1}{T - 1}\right) \overline{\theta}_\Delta, t=1,\ldots,T\right\}
\end{equation*}
for some integer number of points $T$ as a polyhedral relaxation of cosine. This relaxation is such that the tangent lines intersect the cosine curve at points equidistant on the interval ${[-\overline{\theta}_\Delta, \overline{\theta}_\Delta]}$.

\begin{figure}[H]
    \centering
    \ifdefined\epsswitch
        \includegraphics{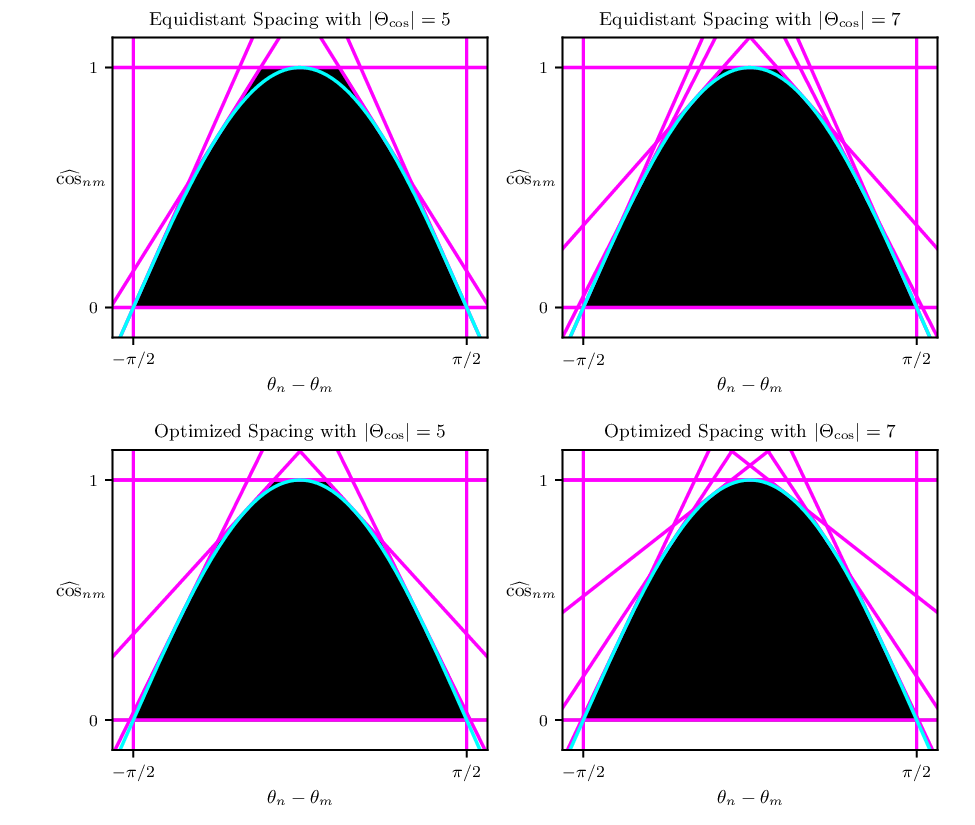}
    \else
        \includegraphics{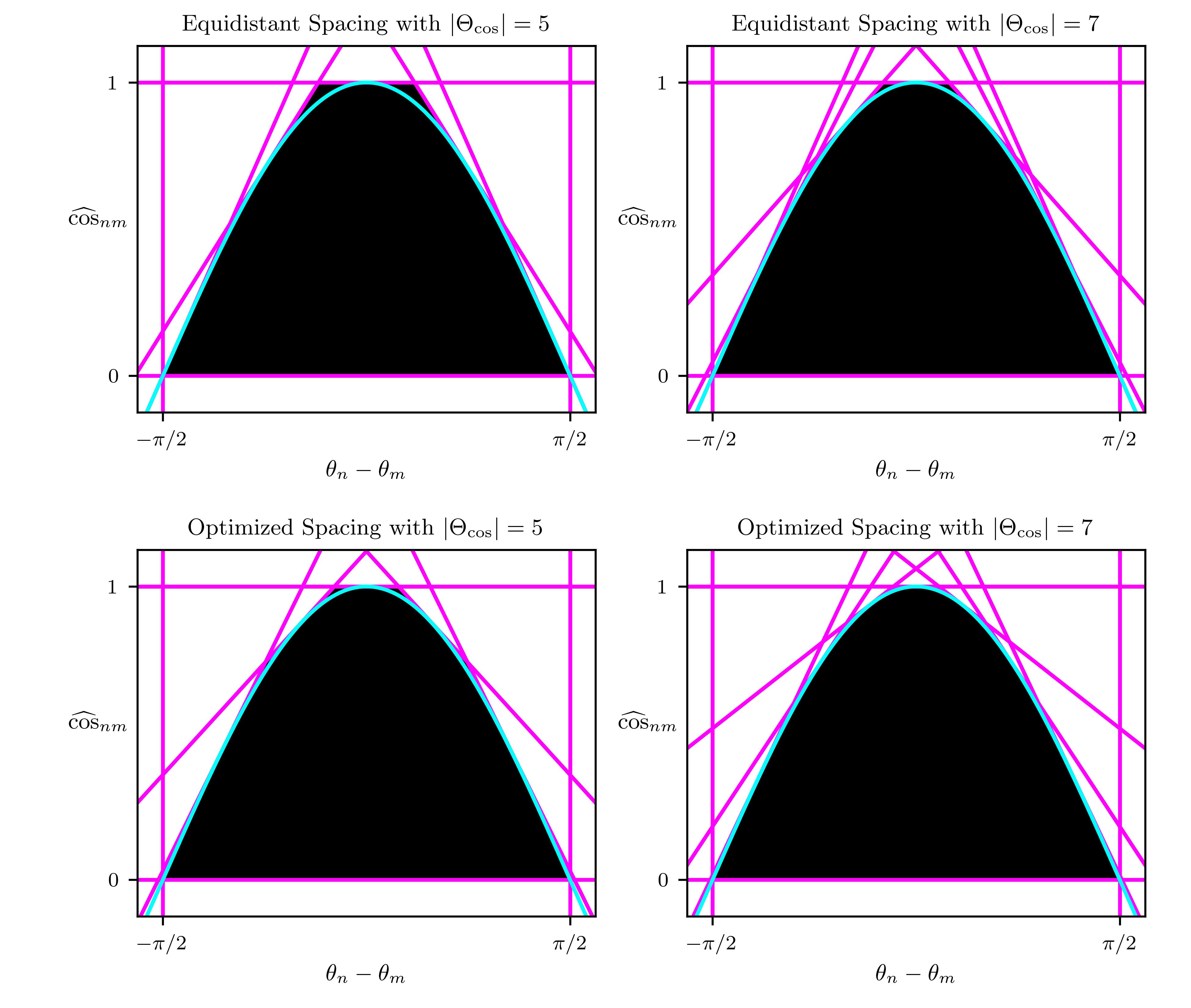}
    \fi
    \caption{
        Comparison of Cosine Relaxations Based on Tangent Lines that are Equidistantly or Optimally Spaced \\[1em]
        \emph{Note: Cyan represents the exact sets, black the additional feasible solutions admitted by the relaxations, and magenta the boundaries of the linear and quadratic constraints that define the relaxed feasible regions.}
    }
    \label{fig:cosine_geometry_comparison}
\end{figure}

A shortcoming of this particular relaxation is that the potential for overestimation error is the greatest in the vicinity of ${\theta = \theta_n - \theta_m = 0}$ as pictured in the top row of subplots in Figure~\ref{fig:cosine_geometry_comparison}. Because our approximation of sine is based on the assumption that ${\theta_n - \theta_m}$ is small, it would be ideal for the relaxation of cosine to be tighter for such values. To remedy this issue, we propose using a $\Theta_\text{cos}$ that minimizes the worst-case relaxation error. Let $\hat{\theta}$ be a point at which some bounding tangent line intersects $\cos(\theta)$ and let ${\Theta_\text{cos} = \{\hat{\theta}_1, \ldots, \hat{\theta}_T\}}$ be a set of $T$ such points.
Subject to ${-\overline{\theta}_\Delta \le \hat{\theta}_1 \le \ldots \le \hat{\theta}_T \le \overline{\theta}_\Delta}$, the set $\Theta_\text{cos}$ that minimizes the maximum relaxation error is the solution to the following trilevel optimization problem:
\begin{equation}
   \min_{\hat{\theta}_1,\ldots,\hat{\theta}_T} \max_{-\overline{\theta}_\Delta \le \theta \le \overline{\theta}_\Delta} \min_{t = 1, \ldots, T}
   \{\mathcal{B}_1(\theta; \hat{\theta}_t) - \cos(\theta)\}.
   \label{eq:tight_cos_trilevel}
\end{equation}
The outer level is a minimum over the decision variables, the middle level is a maximum over the range of permissible $\theta$, and the inner level simply captures the distance between $\cos(\theta)$ and the nearest bounding tangent line at $\theta$.

Intuitively, the worst-case error of a polygonal relaxation may only occur at one of the points where adjacent tangent lines intersect or at one of the endpoints of the interval ${[-\overline{\theta}_\Delta, \overline{\theta}_\Delta]}$. This observation allows us to analytically eliminate the inner-level minimum to simplify the formulation. Let $\tilde{\theta}_{t,t+1}$ be the point $\theta$ at which ${\mathcal{B}_1(\theta; \hat{\theta}_t) = \mathcal{B}_1(\theta; \hat{\theta}_{t+1})}$. Analytically,
\begin{equation}
    \tilde{\theta}_{t,t+1} =
        \frac{\hat{\theta}_t \sin(\hat{\theta}_t) - \hat{\theta}_{t+1} \sin(\hat{\theta}_{t+1}) + \cos(\hat{\theta}_t) - \cos(\hat{\theta}_{t+1})}
             {\sin(\hat{\theta}_t) - \sin(\hat{\theta}_{t+1})}, \label{eq:theta_intersection}
\end{equation}
and
\begin{equation}
    \mathcal{B}_1(\tilde{\theta}_{t,t+1}; \hat{\theta}_t) =
        \frac{\sin(\hat{\theta}_{t}) \cos(\hat{\theta}_{t+1}) - \cos(\hat{\theta}_t) \sin(\hat{\theta}_{t+1}) - (\hat{\theta}_t - \hat{\theta}_{t+1}) \sin(\hat{\theta}_t) \sin(\hat{\theta}_{t+1})}
             {\sin(\hat{\theta}_t) - \sin(\hat{\theta}_{t+1})}. \label{eq:f_theta_intersection}
\end{equation}
This leads to the following reformulation of \eqref{eq:tight_cos_trilevel}:
\begin{mini*}[2]
    {}{z}{}{}
    \addConstraint{-\overline{\theta}_\Delta \le \hat{\theta}_1 \le \ldots \le \hat{\theta}_T \le \overline{\theta}_\Delta}
    \addConstraint{z \ge \mathcal{B}_1(\tilde{\theta}_{t,t+1}; \hat{\theta}_t) - \cos(\tilde{\theta}_{t,t+1}),}
                  {\quad}{t=1,\ldots,T-1,}
    \addConstraint{z \ge \mathcal{B}_1(-\overline{\theta}_\Delta; \hat{\theta}_1) - \cos(-\overline{\theta}_\Delta),}
    \addConstraint{z \ge \mathcal{B}_1(\overline{\theta}_\Delta; \hat{\theta}_T) - \cos(\overline{\theta}_\Delta),}
    \addConstraint{\eqref{eq:theta_intersection}, \eqref{eq:f_theta_intersection}} \notag
\end{mini*}
Here, the second-level maximum is modeled using the variable $z$. This model still involves many trigonometric functions and is accordingly nonconvex. We solve the model for ${|\Theta_\text{cos}| = 5}$ and ${|\Theta_\text{cos}| = 7}$ and for ${\overline{\theta}_\Delta = \tfrac{\pi}{2}}$ using the Knitro solver \citep{Byrd2006} available through the Network-Enabled Optimization Solver (NEOS) Server \citep{Czyzyk1998,Dolan2001,Gropp1997}. For these parameters, the differences between $\Theta_\text{cos}$ being composed of equidistant points versus points spaced as prescribed by our optimization model are illustrated in Figure~\ref{fig:cosine_geometry_comparison}. We use the set comprising 7 points, ${\Theta_\text{cos} = \left\{0, \pm 0.354, \pm 0.735, \pm 1.211\right\}}$, in the implementation of our LPAC-F model.

\section{Tightening the Ohm's Law Constraints}

In our adapted DC and LPAC models, out-of-service branches are treated as open circuits, and we use the big-$M$ modeling technique to only enforce Ohm's Law for operational branches. Unfortunately, this technique requires careful calibration of the big-$M$ constants. Having big-$M$ constants that are too large may cause the solution space of the linear programming (LP) relaxation to be loose relative to the convex hull of the integer-feasible solution space. This affects the bound provided by the LP relaxation subproblems in mixed-integer algorithms like branch-and-bound. In power systems literature, obtaining the tightest big-$M$ value requires solving an optimization problem that incorporates all of the power flow constraints except the Ohm's Law constraint for which the big-$M$ value is being determined, which is rather considered in the objective function. A variety of methods have been proposed for solving these so-called bounding problems \citep{Fattahi2019,Porras2023,Pineda2023}.

For our adaptations of the DC and LPAC models, we solve significantly relaxed versions of these bounding problems. We incorporate the same objective function (\textit{i.e.}, the right-hand side of each Ohm's Law constraint) but only consider a subset of the other constraints and variables in the power flow problem. Specifically, for branch $l \in L_{nm}$, we include only the infeasibility indicator variable $\chi$ and the variables that pertain directly to branch $l$ or buses $n$ and $m$. Similarly, we include only the constraints involving exclusively those variables. We fix $\beta_{nm} = 0$ in these problems since the big-$M$ constant is only applicable in such a case. For our adapted DC model, the optimization problems we use to determine the big-$M$ constants are
\begin{mini*}[2]
    {}{-b_l \widehat{\sin}_{nm}}{}{M_{nml}^\text{active} =}
    \addConstraint{\widehat{\sin} \in \mathcal{Y}_\text{sin}(\theta_n, \theta_m, 0)},
\end{mini*}
for all $(n,m) \in E, l \in L_{nm}$. For our adapted LPAC model, we solve
\begin{mini*}[2]
	{}{v_n g_l \left(v_m - v_n\right) \chi + v_n^2 g_l - v_n v_m \left(g_l \widehat{\text{cos}}_{nm} + b_l \widehat{\text{sin}}_{nm}\right)}{}{M_{nml}^\text{active} =}
    \addConstraint{\widehat{\text{sin}}_{nm} \in \mathcal{Y}_\text{sin}(\theta_n, \theta_m, 0),~
                   \widehat{\text{cos}}_{nm} \in \mathcal{Y}_\text{cos}(\theta_n, \theta_m, 0),}
\end{mini*}
\begin{mini*}[2]
	{}{v_n b_l (v_n - v_m) \chi - v_n^2 b_l - v_n v_m \left(g_l \widehat{\text{sin}}_{nm} - b_l \widehat{\text{cos}}_{nm}\right)}{}{M_{nml}^\text{reactive} =}
    \breakObjective{\phantom{M_{nml}^\text{reactive} =~}
                    \quad - v_n b_l \left(\phi_n - \phi_m\right) - (v_n - v_m) b_l \phi_n}
    \addConstraint{\widehat{\text{sin}}_{nm} \in \mathcal{Y}_\text{sin}(\theta_n, \theta_m, 0),~
                   \widehat{\text{cos}}_{nm} \in \mathcal{Y}_\text{cos}(\theta_n, \theta_m, 0),}
    \addConstraint{\underline{v}_n \le v_n + \phi_n \le \overline{v}_n,~
                   \underline{v}_m \le v_m + \phi_m \le \overline{v}_m.}
\end{mini*}
In these bounding problems, variables are also individually restricted to the domains specified in Section 3.1 of the article. The superscript on each big-$M$ indicates the constraint for which the constant is applicable. As formulated, these minimization problems determine a big-$M$ constant that is valid for the lower bound in the corresponding constraint. To obtain a big-$M$ constant that is valid for the upper bound, the bounding problem may be solved as a maximization problem. The big-$M$ constants are consequential only if $\beta_{nm} = 0$, and the constraints from $\mathcal{L_\text{DC}}$ and $\mathcal{L_\text{LPAC}}$ that are present in the bounding problems simplify in this case. If $\beta_{nm} = 0$, then $\tilde{p}_l = \tilde{q}_l = 0$ hence the absence of $\tilde{p}_l$ and $\tilde{q}_l$ in the bounding problems. Also, if $\beta_{nm} = 0$, the approximation of sine in our adapted DC and LPAC models allows $\widehat{\sin}_{nm} \in [-2 \overline{\theta}, 2 \overline{\theta}]$, the approximation of cosine in our adapted LPAC-C model requires $\widehat{\cos}_{nm} = 1$, and the relaxation of cosine in our adapted LPAC-F and QPAC models allows ${\widehat{\cos}_{nm} \in [\cos(\overline{\theta}_\Delta), 1]}$. That is, $\theta_n$ and $\theta_m$ are eliminated in the simplified bounding problems, and $\widehat{\sin}_{nm}$ and $\widehat{\cos}_{nm}$ are independently constrained. The other variables in these linear big-$M$ optimization problems (\textit{i.e.}, $\phi_n$, $\phi_m$, and $\chi$) are also independently constrained. While including more of the problem context in these optimization problems would produce superior big-$M$ constants, such problems are more difficult to solve. In contrast, our problems are easily solved by decomposition (\textit{i.e.}, by optimizing each of $\widehat{\sin}_{nm}, \widehat{\cos}_{nm}, \phi_n, \phi_m$, and $\chi$ independently).
\section{Flooding Scenarios}

\begin{figure}[H]
    \centering
    \includegraphics{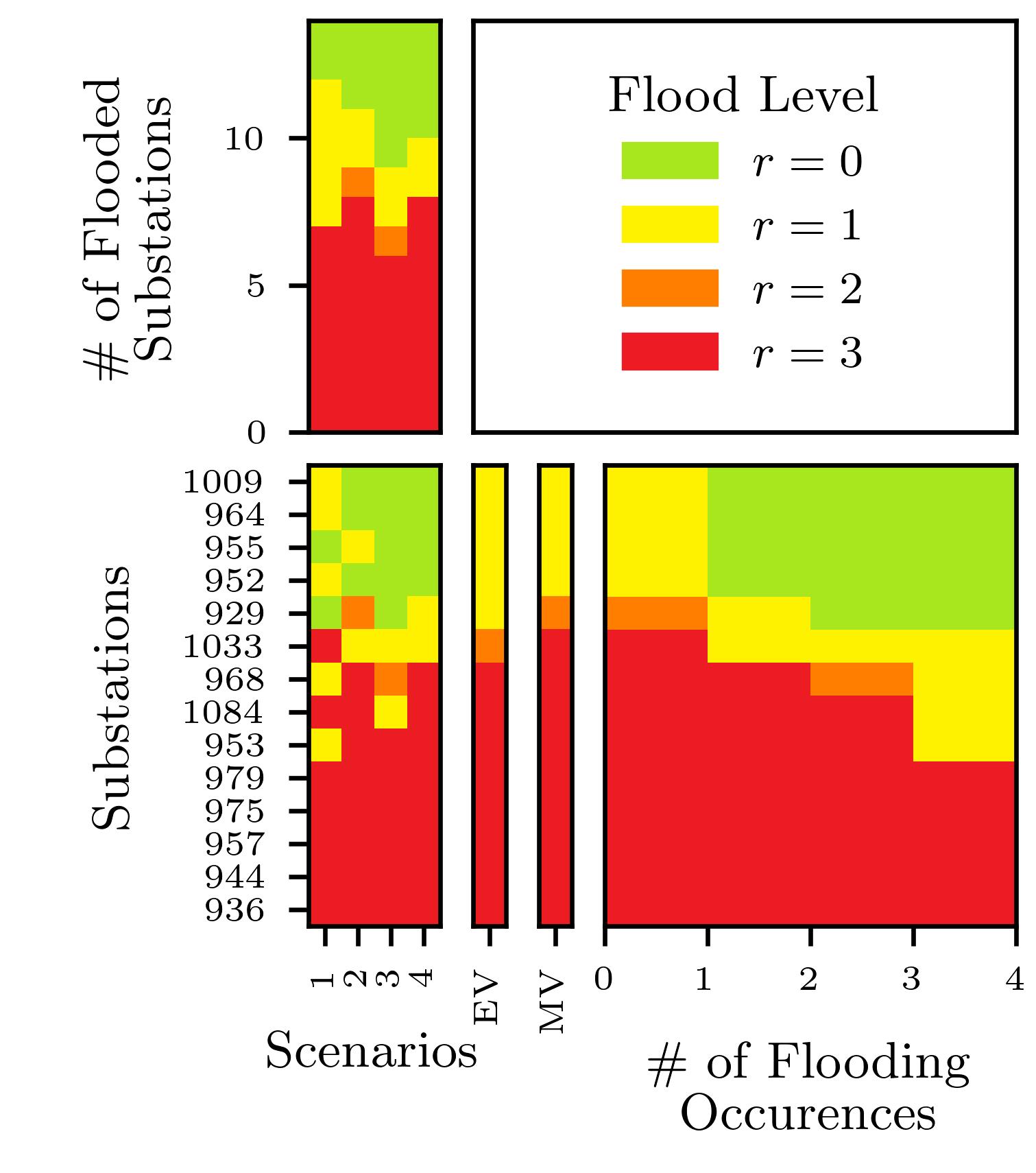}
    \caption{Tropical Storm Imelda Flood Levels by Scenario and by Substation}
    \label{fig:imelda_uncertainty_heatmap}
\end{figure}

\begin{figure}[H]
    \centering
    \includegraphics{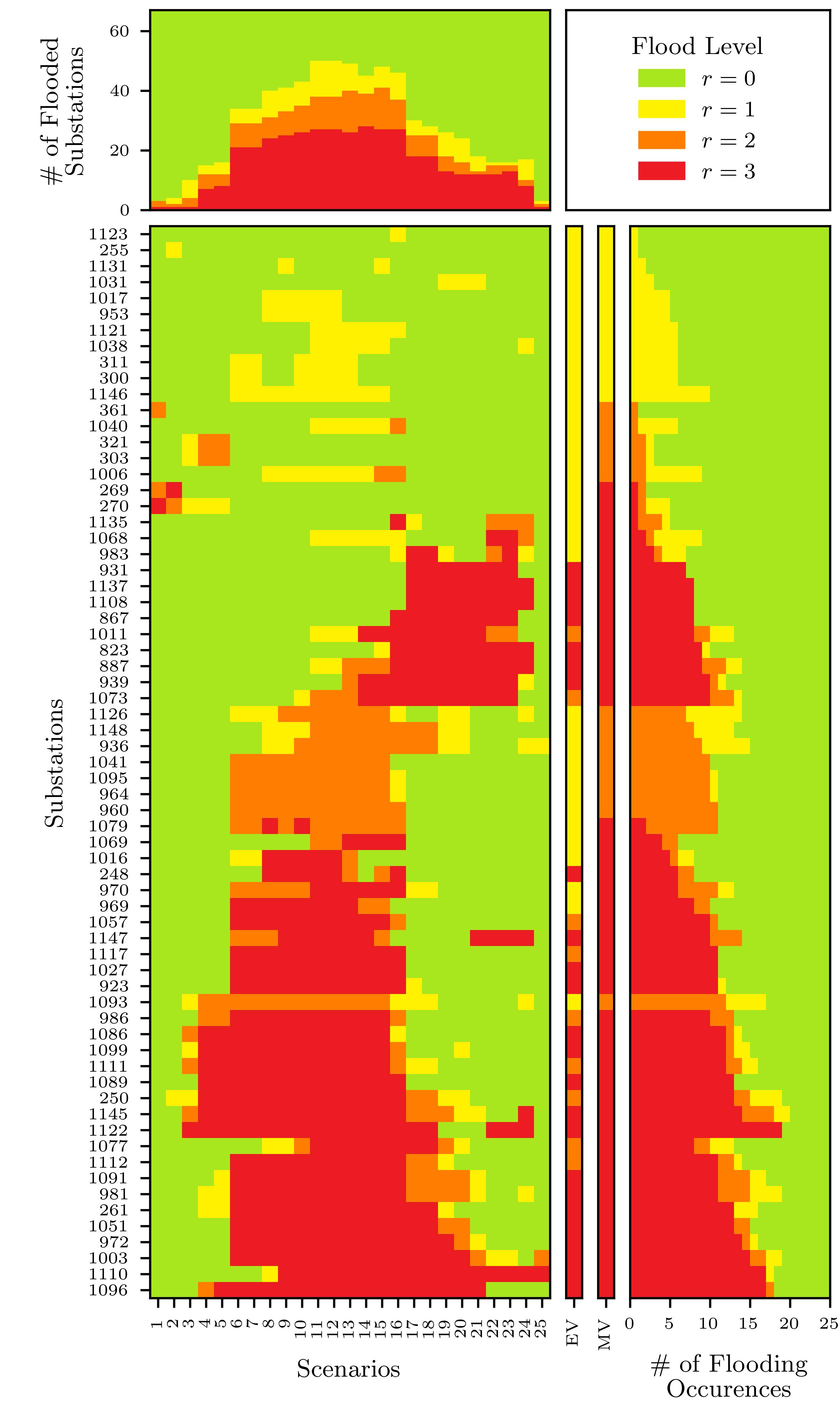}
    \caption{Hurricane Harvey Flood Levels by Scenario and by Substation}
    \label{fig:harvey_uncertainty_heatmap}
\end{figure}

In Figures~\ref{fig:imelda_uncertainty_heatmap} and \ref{fig:harvey_uncertainty_heatmap}, the subplot in the bottom left illustrates the flooding at each affected substation in each scenario. The two subplots to its right show the same for the expected value (EV) and maximum value (MV) scenarios. The top-most subplot is a vertically stacked bar chart indicating how many substations experience each level of flooding in each scenario, and the right-most plot is a horizontally stacked bar chart indicating the number of scenarios in which each substation experiences each flood level.
\section{Tabulated Differences in Solutions}

In Section 4.4.1 of the article, we discuss the impact of the chosen power flow model on the optimal mitigation solutions. Among the 429 comparisons we conducted, there were only 29 cases in which the solutions differed. Those differences are detailed in Table~\ref{tab:solution_differences}. In the table, the two columns with relative optimality gap values are for the solution induced by the former power flow model evaluated in the latter (left) and vice versa (right).

\begin{table}[H]
    \centering
    \small
    \caption{Instances for Which the Mitigation Solutions Differed in the Incorporated Power Flow Model}
    \label{tab:solution_differences}
    \vspace{0.5em}
    \begin{tabular}{|c|c|r|C{1.55cm}|C{1.55cm}|R{1.45cm}|R{1.45cm}|R{1.15cm}|R{1.15cm}|}
        \hline
        Case Study & Model & Budget & \multicolumn{2}{c|}{Power Flow Models} & AbsSim & RelSim & \multicolumn{2}{c|}{Opt. Gaps (\%)} \\
        \hline
        Imelda & SP &   4 & DC     & LPAC-F &   2 & 0.5000 & 0.9730 & 1.1138 \\
        Imelda & SP &   4 & DC     & QPAC   &   2 & 0.5000 & 0.1559 & 1.1138 \\
        Imelda & SP &   4 & LPAC-C & LPAC-F &   2 & 0.5000 & 0.9730 & 1.1284 \\
        Imelda & SP &   4 & LPAC-C & QPAC   &   2 & 0.5000 & 0.1559 & 1.1284 \\
        \hline
        Imelda & RO &   2 & DC     & LPAC-F &   1 & 0.5000 & 2.6918 & 3.6036 \\
        Imelda & RO &   2 & DC     & QPAC   &   1 & 0.5000 & 0.2754 & 3.6036 \\
        Imelda & RO &   2 & LPAC-C & LPAC-F &   1 & 0.5000 & 2.6918 & 3.5914 \\
        Imelda & RO &   2 & LPAC-C & QPAC   &   1 & 0.5000 & 0.2754 & 3.5914 \\
        Imelda & RO &   2 & LPAC-F & QPAC   &   1 & 0.5000 & 0.0011 & 0.0000 \\
        \hline
        Harvey & SP &  88 & DC     & LPAC-C &  85 & 0.9659 & 0.1100 & 0.1082 \\
        Harvey & SP & 110 & DC     & LPAC-C & 107 & 0.9727 & 0.0175 & 0.0179 \\
        Harvey & SP & 114 & DC     & LPAC-C & 110 & 0.9649 & 0.0009 & 0.0008 \\
        Harvey & SP & 117 & DC     & LPAC-C & 113 & 0.9658 & 0.0009 & 0.0008 \\
        Harvey & SP & 128 & DC     & LPAC-C & 118 & 0.9219 & 0.0137 & 0.0049 \\
        Harvey & SP & 130 & DC     & LPAC-C & 127 & 0.9769 & 0.0006 & 0.0009 \\
        Harvey & SP & 137 & DC     & LPAC-C & 133 & 0.9708 & 0.0061 & 0.0026 \\
        Harvey & SP & 140 & DC     & LPAC-C & 137 & 0.9785 & 0.0064 & 0.0025 \\
        Harvey & SP & 144 & DC     & LPAC-C & 141 & 0.9792 & 0.0064 & 0.0025 \\
        Harvey & SP & 148 & DC     & LPAC-C & 141 & 0.9527 & 0.0031 & 0.0056 \\
        Harvey & SP & 150 & DC     & LPAC-C & 147 & 0.9800 & 0.0064 & 0.0024 \\
        Harvey & SP & 156 & DC     & LPAC-C & 153 & 0.9808 & 0.0005 & 0.0007 \\
        Harvey & SP & 164 & DC     & LPAC-C & 161 & 0.9817 & 0.0021 & 0.0044 \\
        Harvey & SP & 177 & DC     & LPAC-C & 174 & 0.9831 & 0.0020 & 0.0020 \\
        Harvey & SP & 186 & DC     & LPAC-C & 183 & 0.9839 & 0.0021 & 0.0200 \\
        \hline
        Harvey & RO &   1 & DC     & LPAC-C &   1 & 0.0000 & 0.0087 & 0.0000 \\
        Harvey & RO &   5 & DC     & LPAC-C &   5 & 0.6000 & 0.0158 & 0.0005 \\
        Harvey & RO &  34 & DC     & LPAC-C &  34 & 0.9118 & 0.0003 & 0.0000 \\
        Harvey & RO &  35 & DC     & LPAC-C &  35 & 0.9143 & 0.0003 & 0.0000 \\
        Harvey & RO &  36 & DC     & LPAC-C &  36 & 0.8889 & 0.0110 & 0.0020 \\
        \hline
    \end{tabular}
\end{table}

\end{appendices}

% Acknowledgments here
\ACKNOWLEDGMENT{%
We the authors thank The University of Texas Energy Institute for funding our ``Defending the Electricity Infrastructure Against Extreme Weather Events, Now \& in the Future'' research project as part of the ``Fueling a Sustainable Energy Transition'' initiative.

We also acknowledge the Texas Advanced Computing Center (TACC) at The University of Texas at Austin for providing HPC and data visualization resources that have contributed to the research results reported within this paper. URL: http://www.tacc.utexas.edu

Finally, we are grateful for several individual contributions to this work by current and former members of our research team at The University of Texas at Austin. Namely, we thank Dr. Zong-Liang Yang, Dr. Carey King, Dr. Wen-Ying Wu, Dr. Kyoung Yoon Kim, Joshua Yip, and Ashutosh Shukla. We also thank the three anonymous reviewers of our manuscript whose detailed comments strengthened the paper significantly.
}% Leave this (end of acknowledgment)

% CASE 1: BiBTeX used to constantly update the references 
%   (while the paper is being written).
\bibliographystyle{informs2014} % outcomment this and next line in Case 1
\bibliography{main} % if more than one, comma separated

\begin{thebibliography}{61}
\providecommand{\natexlab}[1]{#1}
\providecommand{\url}[1]{\texttt{#1}}
\providecommand{\urlprefix}{URL }

\bibitem[{Arab et~al.(2015)Arab, Khodaei, Khator, Ding, Emesih, \protect\BIBand{} Han}]{Arab2015}
Arab A, Khodaei A, Khator SK, Ding K, Emesih VA, Han Z (2015) Stochastic pre-hurricane restoration planning for electric power systems infrastructure. \emph{IEEE Transactions on Smart Grid} 6(2):1046--1054, \urlprefix\url{http://dx.doi.org/10.1109/TSG.2015.2388736}.

\bibitem[{Asghari et~al.(2022)Asghari, Fathollahi-Fard, Mirzapour Al-e hashem, \protect\BIBand{} Dulebenets}]{Asghari2022}
Asghari M, Fathollahi-Fard AM, Mirzapour Al-e hashem SMJ, Dulebenets MA (2022) Transformation and linearization techniques in optimization: A state-of-the-art survey. \emph{Mathematics} 10(2), \urlprefix\url{http://dx.doi.org/10.3390/math10020283}.

\bibitem[{Austgen et~al.(2021)Austgen, Hasenbein, \protect\BIBand{} Kutanoglu}]{Austgen2021}
Austgen B, Hasenbein J, Kutanoglu E (2021) Impacts of approximate power flow models on optimal flood mitigation in a stochastic program. \emph{{IIE} {Annual} {Conference} {Proceedings}}, 518--523.

\bibitem[{Austgen et~al.(2022)Austgen, Kutanoglu, \protect\BIBand{} Hasenbein}]{Austgen2022a}
Austgen B, Kutanoglu E, Hasenbein JJ (2022) Comparison of linear power flow approximations in the context of winter storm planning. \emph{{IIE} {Annual} {Conference} {Proceedings}}, 1--6.

\bibitem[{Austgen et~al.(2023)Austgen, Kutanoglu, \protect\BIBand{} Hasenbein}]{Austgen2023}
Austgen B, Kutanoglu E, Hasenbein JJ (2023) A two-stage stochastic programming model for electric substation flood mitigation prior to an imminent hurricane, \urlprefix\url{http://dx.doi.org/10.48550/ARXIV.2302.10996}, {Submitted}.

\bibitem[{Austgen et~al.(2024)Austgen, Kutanoglu, Hasenbein, \protect\BIBand{} Santoso}]{Austgen2024code}
Austgen B, Kutanoglu E, Hasenbein JJ, Santoso S (2024) Comparisons of two-stage models for flood mitigation of electrical substations. \urlprefix\url{http://dx.doi.org/10.1287/ijoc.2023.0125.cd}, available for download at https://github.com/INFORMSJoC/2023.0125.

\bibitem[{Ben-Tal et~al.(2009)Ben-Tal, Ghaoui, \protect\BIBand{} Nemirovski}]{BenTal2009}
Ben-Tal A, Ghaoui L, Nemirovski A (2009) \emph{Robust {Optimization}}. Princeton {Series} in {Applied} {Mathematics} (Princeton University Press).

\bibitem[{Birge \protect\BIBand{} Louveaux(2011)}]{Birge2011}
Birge JR, Louveaux F (2011) \emph{Introduction to {Stochastic} {Programming}} (Springer Science \& Business Media).

\bibitem[{Blake \protect\BIBand{} Zelinsky(2018)}]{Blake2018}
Blake ES, Zelinsky DA (2018) Hurricane {Harvey}. Technical Report AL092017, National Hurricane Center, \urlprefix\url{https://www.nhc.noaa.gov/data/tcr/AL092017_Harvey.pdf}.

\bibitem[{Byrd et~al.(2006)Byrd, Nocedal, \protect\BIBand{} Waltz}]{Byrd2006}
Byrd RH, Nocedal J, Waltz RA (2006) Knitro: {An} integrated package for nonlinear optimization. Di~Pillo G, Roma M, eds., \emph{Large-{Scale} {Nonlinear} {Optimization}}, 35--59 (Boston, MA: Springer US), \urlprefix\url{http://dx.doi.org/10.1007/0-387-30065-1_4}.

\bibitem[{Coffrin et~al.(2019)Coffrin, Bent, Tasseff, Sundar, \protect\BIBand{} Backhaus}]{Coffrin2019}
Coffrin C, Bent R, Tasseff B, Sundar K, Backhaus S (2019) Relaxations of ac maximal load delivery for severe contingency analysis. \emph{IEEE Transactions on Power Systems} 34(2):1450--1458, \urlprefix\url{http://dx.doi.org/10.1109/TPWRS.2018.2876507}.

\bibitem[{Coffrin et~al.(2011)Coffrin, Hentenryck, \protect\BIBand{} Bent}]{Coffrin2011}
Coffrin C, Hentenryck PV, Bent R (2011) Strategic stockpiling of power system supplies for disaster recovery. \emph{2011 {IEEE} {Power} and {Energy} {Society} {General} {Meeting}}, 1--8.

\bibitem[{Coffrin et~al.(2020)Coffrin, Hijazi, \protect\BIBand{} Van~Hentenryck}]{Coffrin2020}
Coffrin C, Hijazi H, Van~Hentenryck P (2020) Alternating current {(AC)} power flow analysis in an electrical power network. US Patent 10,591,520.

\bibitem[{Coffrin \protect\BIBand{} Van~Hentenryck(2014)}]{Coffrin2014}
Coffrin C, Van~Hentenryck P (2014) A linear-programming approximation of {AC} power flows. \emph{INFORMS Journal on Computing} 26(4):718--734, \urlprefix\url{http://dx.doi.org/10.1287/ijoc.2014.0594}.

\bibitem[{Coffrin \protect\BIBand{} Van~Hentenryck(2015)}]{Coffrin2015}
Coffrin C, Van~Hentenryck P (2015) Transmission system restoration with co-optimization of repairs, load pickups, and generation dispatch. \emph{International Journal of Electrical Power \& Energy Systems} 72:144--154, \urlprefix\url{http://dx.doi.org/10.1016/j.ijepes.2015.02.027}.

\bibitem[{{Czyzyk} et~al.(1998){Czyzyk}, {Mesnier}, \protect\BIBand{} {Mor\'e}}]{Czyzyk1998}
{Czyzyk} J, {Mesnier} MP, {Mor\'e} JJ (1998) The {NEOS} server. \emph{IEEE Journal on Computational Science and Engineering} 5(3):68--75.

\bibitem[{{Dolan}(2001)}]{Dolan2001}
{Dolan} ED (2001) The {NEOS} server 4.0 administrative guide. Technical Memorandum ANL/MCS-TM-250, Mathematics and Computer Science Division, Argonne National Laboratory.

\bibitem[{Dvorkin et~al.(2018)Dvorkin, Henneaux, Kirschen, \protect\BIBand{} Pandžić}]{Dvorkin2018}
Dvorkin Y, Henneaux P, Kirschen DS, Pandžić H (2018) Optimizing primary response in preventive security-constrained optimal power flow. \emph{IEEE Systems Journal} 12(1):414--423, \urlprefix\url{http://dx.doi.org/10.1109/JSYST.2016.2527726}.

\bibitem[{Eijgenraam et~al.(2017)Eijgenraam, Brekelmans, den Hertog, \protect\BIBand{} Roos}]{Eijgenraam2017}
Eijgenraam C, Brekelmans R, den Hertog D, Roos K (2017) Optimal strategies for flood prevention. \emph{Management Science} 63(5):1644--1656, \urlprefix\url{http://dx.doi.org/10.1287/mnsc.2015.2395}.

\bibitem[{Eijgenraam et~al.(2014)Eijgenraam, Kind, Bak, Brekelmans, den Hertog, Duits, Roos, Vermeer, \protect\BIBand{} Kuijken}]{Eijgenraam2014}
Eijgenraam C, Kind J, Bak C, Brekelmans R, den Hertog D, Duits M, Roos K, Vermeer P, Kuijken W (2014) Economically efficient standards to protect the netherlands against flooding. \emph{Interfaces} 44(1):7--21, \urlprefix\url{http://dx.doi.org/10.1287/inte.2013.0721}.

\bibitem[{Fattahi et~al.(2019)Fattahi, Lavaei, \protect\BIBand{} Atamtürk}]{Fattahi2019}
Fattahi S, Lavaei J, Atamtürk A (2019) A bound strengthening method for optimal transmission switching in power systems. \emph{IEEE Transactions on Power Systems} 34(1):280--291, \urlprefix\url{http://dx.doi.org/10.1109/TPWRS.2018.2867999}.

\bibitem[{Garcia et~al.(2022)Garcia, Austgen, Pierre, Hasenbein, \protect\BIBand{} Kutanoglu}]{Garcia2022}
Garcia M, Austgen B, Pierre B, Hasenbein J, Kutanoglu E (2022) Risk-averse investment optimization for power system resilience to winter storms. \emph{{IEEE}/{PES} {Transmission} and {Distribution} {Conference} and {Exposition}}.

\bibitem[{Garifi et~al.(2022)Garifi, Johnson, Arguello, \protect\BIBand{} Pierre}]{Garifi2022}
Garifi K, Johnson ES, Arguello B, Pierre BJ (2022) Transmission grid resiliency investment optimization model with {SOCP} recovery planning. \emph{IEEE Transactions on Power Systems} 37(1):26--37, \urlprefix\url{http://dx.doi.org/10.1109/TPWRS.2021.3091538}.

\bibitem[{{Gropp} \protect\BIBand{} {Mor\'e}(1997)}]{Gropp1997}
{Gropp} W, {Mor\'e} JJ (1997) Optimization environments and the {NEOS} server. {Buhman} MD, {Iserles} A, eds., \emph{Approximation Theory and Optimization}, 167 -- 182 (Cambridge University Press).

\bibitem[{{Gurobi Optimization, LLC}(2022)}]{Gurobi2022}
{Gurobi Optimization, LLC} (2022) Gurobi {Optimizer} {Reference} {Manual}. \urlprefix\url{https://www.gurobi.com}.

\bibitem[{Huang et~al.(2022)Huang, Lai, Zhao, Yang, Zhong, \protect\BIBand{} Lai}]{Huang2022}
Huang L, Lai CS, Zhao Z, Yang G, Zhong B, Lai LL (2022) Robust {N}-k security-constrained optimal power flow incorporating preventive and corrective generation dispatch to improve power system reliability. \emph{CSEE Journal of Power and Energy Systems} 1--14, \urlprefix\url{http://dx.doi.org/10.17775/CSEEJPES.2021.06560}.

\bibitem[{IEEE2022()}]{IEEEReliabilityIndices2022}
IEEE2022 (2022) {IEEE} guide for electric power distribution reliability indices. \emph{IEEE Std 1366-2022 (Revision of IEEE Std 1366-2012)} 1--44, \urlprefix\url{http://dx.doi.org/10.1109/IEEESTD.2022.9955492}.

\bibitem[{ISER()}]{ISER2023}
ISER (2022) {ISER} {Electric} {Disturbance} {Events} ({DOE}-417). \urlprefix\url{www.oe.netl.doe.gov/oe417.aspx}.

\bibitem[{Kile et~al.(2014)Kile, Uhlen, Warland, \protect\BIBand{} Kjølle}]{Kile2014}
Kile H, Uhlen K, Warland L, Kjølle G (2014) A comparison of {AC} and {DC} power flow models for contingency and reliability analysis. \emph{2014 {Power} {Systems} {Computation} {Conference}}, 1--7, \urlprefix\url{http://dx.doi.org/10.1109/PSCC.2014.7038459}.

\bibitem[{Klerk et~al.(2021)Klerk, Kanning, Kok, \protect\BIBand{} Wolfert}]{Klerk2021}
Klerk WJ, Kanning W, Kok M, Wolfert R (2021) Optimal planning of flood defence system reinforcements using a greedy search algorithm. \emph{Reliability Engineering \& System Safety} 207:107344, ISSN 0951-8320, \urlprefix\url{http://dx.doi.org/https://doi.org/10.1016/j.ress.2020.107344}.

\bibitem[{Knueven et~al.(2020)Knueven, Ostrowski, \protect\BIBand{} Watson}]{Knueven2020}
Knueven B, Ostrowski J, Watson JP (2020) On mixed-integer programming formulations for the unit commitment problem. \emph{INFORMS Journal on Computing} 32(4):857--876, \urlprefix\url{http://dx.doi.org/10.1287/ijoc.2019.0944}.

\bibitem[{Knutson et~al.(2020)Knutson, Camargo, Chan, Emanuel, Ho, Kossin, Mohapatra, Satoh, Sugi, Walsh, \protect\BIBand{} Wu}]{Knutson2020}
Knutson T, Camargo SJ, Chan JCL, Emanuel K, Ho CH, Kossin J, Mohapatra M, Satoh M, Sugi M, Walsh K, Wu L (2020) Tropical cyclones and climate change assessment: {Part} {II}: Projected response to anthropogenic warming. \emph{Bulletin of the American Meteorological Society} 101(3):E303 -- E322, \urlprefix\url{http://dx.doi.org/10.1175/BAMS-D-18-0194.1}.

\bibitem[{Latto \protect\BIBand{} Berg(2020)}]{Latto2020}
Latto A, Berg R (2020) Tropical {Storm} {Imelda}. Technical Report AL112019, National Hurricane Center, \urlprefix\url{https://www.nhc.noaa.gov/data/tcr/AL112019_Imelda.pdf}.

\bibitem[{Linkov \protect\BIBand{} Trump(2019)}]{Linkov2019}
Linkov I, Trump BD (2019) \emph{The Science and Practice of Resilience} (Springer), \urlprefix\url{https://link.springer.com/book/10.1007/978-3-030-04565-4}.

\bibitem[{Liu et~al.(2009)Liu, Tesfatsion, \protect\BIBand{} Chowdhury}]{Liu2009}
Liu H, Tesfatsion L, Chowdhury AA (2009) Locational marginal pricing basics for restructured wholesale power markets. \emph{2009 IEEE Power \& Energy Society General Meeting}, 1--8, \urlprefix\url{http://dx.doi.org/10.1109/PES.2009.5275503}.

\bibitem[{Logan et~al.(2022)Logan, Aven, Guikema, \protect\BIBand{} Flage}]{Logan2022}
Logan TM, Aven T, Guikema SD, Flage R (2022) Risk science offers an integrated approach to resilience. \emph{Nature Sustainability} 5(9):741--748, \urlprefix\url{http://dx.doi.org/10.1038/s41893-022-00893-w}.

\bibitem[{Mohagheghi \protect\BIBand{} Rebennack(2015)}]{Mohagheghi2015}
Mohagheghi S, Rebennack S (2015) Optimal resilient power grid operation during the course of a progressing wildfire. \emph{International Journal of Electrical Power \& Energy Systems} 73:843--852, \urlprefix\url{http://dx.doi.org/10.1016/j.ijepes.2015.05.035}.

\bibitem[{Molzahn \protect\BIBand{} Hiskens(2019)}]{Molzahn2019}
Molzahn DK, Hiskens IA (2019) A {Survey} of {Relaxations} and {Approximations} of the {Power} {Flow} {Equations}. \emph{Foundations and Trends in Electric Energy Systems} 4(1-2):1--221, \urlprefix\url{http://dx.doi.org/10.1561/3100000012}.

\bibitem[{Movahednia \protect\BIBand{} Kargarian(2022)}]{Movahednia2022b}
Movahednia M, Kargarian A (2022) Flood-aware optimal power flow for proactive day-ahead transmission substation hardening. \emph{2022 {IEEE} {Texas} {Power} and {Energy} {Conference} ({TPEC})}, 1--5, \urlprefix\url{http://dx.doi.org/10.1109/TPEC54980.2022.9750830}.

\bibitem[{Movahednia et~al.(2022)Movahednia, Kargarian, Ozdemir, \protect\BIBand{} Hagen}]{Movahednia2022a}
Movahednia M, Kargarian A, Ozdemir CE, Hagen SC (2022) Power grid resilience enhancement via protecting electrical substations against flood hazards: A stochastic framework. \emph{IEEE Transactions on Industrial Informatics} 18(3):2132--2143, \urlprefix\url{http://dx.doi.org/10.1109/TII.2021.3100079}.

\bibitem[{NCEI2022()}]{BillionDollarDisasters2022}
NCEI2022 (2022) U.{S}. billion-dollar weather and climate disasters. Technical report, NOAA National Centers for Environmental Information (NCEI), \urlprefix\url{http://dx.doi.org/10.25921/stkw-7w73}.

\bibitem[{NRC2012()}]{NRC2012}
NRC2012 (2012) \emph{Disaster resilience: {A} national imperative}. \urlprefix\url{http://dx.doi.org/10.17226/13457}.

\bibitem[{Overbye et~al.(2004)Overbye, Cheng, \protect\BIBand{} Sun}]{Overbye2004}
Overbye T, Cheng X, Sun Y (2004) A comparison of the {AC} and {DC} power flow models for {LMP} calculations. \emph{37th {Annual} {Hawaii} {International} {Conference} on {System} {Sciences}, 2004. {Proceedings} of the}, \urlprefix\url{http://dx.doi.org/10.1109/HICSS.2004.1265164}.

\bibitem[{Paul et~al.(2017)Paul, Pathak, Pal, \protect\BIBand{} Chanda}]{Paul2017}
Paul B, Pathak MK, Pal J, Chanda CK (2017) A comparison of locational marginal prices and locational load shedding marginal prices in a deregulated competitive power market. \emph{2017 IEEE Calcutta Conference (CALCON)}, 46--50, \urlprefix\url{http://dx.doi.org/10.1109/CALCON.2017.8280693}.

\bibitem[{Pierre et~al.(2018)Pierre, Arguello, Staid, \protect\BIBand{} Guttromson}]{Pierre2018}
Pierre BJ, Arguello B, Staid A, Guttromson RT (2018) Investment optimization to improve power system resilience. \emph{2018 {IEEE} {International} {Conference} on {Probabilistic} {Methods} {Applied} to {Power} {Systems} ({PMAPS})}, 1--6, \urlprefix\url{http://dx.doi.org/10.1109/PMAPS.2018.8440467}.

\bibitem[{Pineda et~al.(2023)Pineda, Morales, Álvaro Porras, \protect\BIBand{} Domínguez}]{Pineda2023}
Pineda S, Morales JM, Álvaro Porras, Domínguez C (2023) Tight big-{M}s for optimal transmission switching.

\bibitem[{Porras et~al.(2023)Porras, Dom\'{\i}nguez, Morales, \protect\BIBand{} Pineda}]{Porras2023}
Porras A, Dom\'{\i}nguez C, Morales JM, Pineda S (2023) Tight and compact sample average approximation for joint chance-constrained problems with applications to optimal power flow. \emph{INFORMS Journal on Computing} \urlprefix\url{http://dx.doi.org/10.1287/ijoc.2022.0302}.

\bibitem[{Quarm et~al.(2022)Quarm, Fan, Elizondo, \protect\BIBand{} Madani}]{Quarm2022}
Quarm E, Fan X, Elizondo M, Madani R (2022) Proactive posturing of large power grid for mitigating hurricane impacts. \emph{2022 {IEEE} {Power} \& {Energy} {Society} {Innovative} {Smart} {Grid} {Technologies} {Conference} ({ISGT})}, 1--5, \urlprefix\url{http://dx.doi.org/10.1109/ISGT50606.2022.9817529}.

\bibitem[{Rhodes et~al.(2021)Rhodes, Fobes, Coffrin, \protect\BIBand{} Roald}]{Rhodes2021}
Rhodes N, Fobes DM, Coffrin C, Roald L (2021) Powermodelsrestoration.jl: An open-source framework for exploring power network restoration algorithms. \emph{Electric Power Systems Research} 190:106736, \urlprefix\url{http://dx.doi.org/https://doi.org/10.1016/j.epsr.2020.106736}.

\bibitem[{Sahraei-Ardakani \protect\BIBand{} Ou(2017)}]{SahraeiArdakani2017}
Sahraei-Ardakani M, Ou G (2017) Day-ahead preventive scheduling of power systems during natuaral hazards via stochastic optimization. \emph{2017 {IEEE} {Power} {Energy} {Society} {General} {Meeting}}, 1--1, \urlprefix\url{http://dx.doi.org/10.1109/PESGM.2017.8274453}.

\bibitem[{Shukla et~al.(2023)Shukla, Kutanoglu, \protect\BIBand{} Hasenbein}]{Shukla2023}
Shukla A, Kutanoglu E, Hasenbein JJ (2023) Scenario-based optimization models for power grid resilience to extreme flooding events.

\bibitem[{Souto et~al.(2022)Souto, Yip, Wu, Austgen, Kutanoglu, Hasenbein, Yang, King, \protect\BIBand{} Santoso}]{Souto2022}
Souto L, Yip J, Wu WY, Austgen B, Kutanoglu E, Hasenbein J, Yang ZL, King CW, Santoso S (2022) Power system resilience to floods: {Modeling}, impact assessment, and mid-term mitigation strategies. \emph{International Journal of Electrical Power \& Energy Systems} 135:107545, \urlprefix\url{http://dx.doi.org/https://doi.org/10.1016/j.ijepes.2021.107545}.

\bibitem[{TACC(2023)}]{Stampede2}
TACC (2023) Stampede2 - {TACC} {HPC} {Documentation}. \urlprefix\url{docs.tacc.utexas.edu/hpc/stampede2/}.

\bibitem[{Tan et~al.(2018)Tan, Das, Arabshahi, \protect\BIBand{} Kirschen}]{Tan2018}
Tan Y, Das AK, Arabshahi P, Kirschen DS (2018) Distribution systems hardening against natural disasters. \emph{IEEE Transactions on Power Systems} 33(6):6849--6860, \urlprefix\url{http://dx.doi.org/10.1109/TPWRS.2018.2836391}.

\bibitem[{Tasseff et~al.(2019)Tasseff, Bent, \protect\BIBand{} Van~Hentenryck}]{Tasseff2019}
Tasseff B, Bent R, Van~Hentenryck P (2019) Optimization of structural flood mitigation strategies. \emph{Water Resources Research} 55(2):1490--1509, \urlprefix\url{http://dx.doi.org/https://doi.org/10.1029/2018WR024362}.

\bibitem[{Wang et~al.(2013)Wang, Watson, \protect\BIBand{} Guan}]{Wang2013}
Wang Q, Watson JP, Guan Y (2013) Two-stage robust optimization for ${N}-k$ contingency-constrained unit commitment. \emph{IEEE Transactions on Power Systems} 28(3):2366--2375, \urlprefix\url{http://dx.doi.org/10.1109/TPWRS.2013.2244619}.

\bibitem[{Watson et~al.(2014)Watson, Guttromson, Silva-Monroy, Jeffers, Jones, Ellison, Rath, Gearhart, Jones, Corbet, Hanley, \protect\BIBand{} Walker}]{Watson2014}
Watson JP, Guttromson R, Silva-Monroy C, Jeffers R, Jones K, Ellison J, Rath C, Gearhart J, Jones D, Corbet T, Hanley C, Walker LT (2014) Conceptual framework for developing resilience metrics for the electricity, oil, and gas sectors in the {United} {States}. Technical report, Sandia National Laboratories, \urlprefix\url{http://dx.doi.org/10.2172/1177743}.

\bibitem[{Webster et~al.(2005)Webster, Holland, Curry, \protect\BIBand{} Chang}]{Webster2005}
Webster PJ, Holland GJ, Curry JA, Chang HR (2005) Changes in tropical cyclone number, duration, and intensity in a warming environment. \emph{Science} 309(5742):1844--1846, \urlprefix\url{http://dx.doi.org/10.1126/science.1116448}.

\bibitem[{Wächter \protect\BIBand{} Biegler(2006)}]{Wachter2006}
Wächter A, Biegler LT (2006) On the implementation of an interior-point filter line-search algorithm for large-scale nonlinear programming. \emph{Mathematical Programming} 106(1):25--57, \urlprefix\url{http://dx.doi.org/10.1007/s10107-004-0559-y}.

\bibitem[{Yang et~al.(2023)Yang, Rhodes, Yang, Roald, \protect\BIBand{} Ntaimo}]{Yang2023}
Yang H, Rhodes N, Yang H, Roald L, Ntaimo L (2023) Multi-period power system risk minimization under wildfire disruptions.

\bibitem[{Zhang et~al.(2012)Zhang, Vittal, Heydt, \protect\BIBand{} Quintero}]{Zhang2012}
Zhang H, Vittal V, Heydt GT, Quintero J (2012) A mixed-integer linear programming approach for multi-stage security-constrained transmission expansion planning. \emph{IEEE Transactions on Power Systems} 27(2):1125--1133, \urlprefix\url{http://dx.doi.org/10.1109/TPWRS.2011.2178000}.

\end{thebibliography}

% CASE 2: BiBTeX used to generate mypaper.bbl (to be further fine tuned)
%\input{mypaper.bbl} % outcomment this line in Case 2

\end{document}